\documentclass{article}
\usepackage{arxiv}

\usepackage[normalem]{ulem} 

\usepackage[utf8]{inputenc} 
\usepackage[T1]{fontenc}    
\usepackage{hyperref}       
\usepackage{url}            
\usepackage{booktabs}       
\usepackage{amsfonts}       
\usepackage{microtype}      
\usepackage{color}
\usepackage{lipsum}
\usepackage{amsthm}
\usepackage{indentfirst}
\usepackage{graphicx}
\usepackage{epstopdf}
\usepackage{fourier}
\usepackage{bm}
\usepackage{bbm} 
\usepackage{mathtools}
\usepackage{latexsym,enumerate}
\usepackage{soul}
\usepackage{subfig}
\usepackage{adjustbox}

\def\bm{\boldsymbol}

\newcommand{\comment}[1]{}

\newcommand{\BEA}{\begin{eqnarray}}
\newcommand{\EEA}{\end{eqnarray}}

\newcommand{\td}{{\rm d}}
\newcommand{\OL}{\mathcal{L}}

\newcommand{\mbv}{\mathbf{v}}

\newcommand{\Ba}{\bm{a}}
\newcommand{\Bb}{\bm{b}}
\newcommand{\Bg}{\bm{g}}
\newcommand{\Bc}{\bm{c}}
\newcommand{\Bh}{\bm{h}}
\newcommand{\Br}{\bm{r}}
\newcommand{\Bw}{\bm{w}}
\newcommand{\Bx}{\bm{x}}
\newcommand{\By}{\bm{y}}

\newcommand{\Bu}{\bm{u}}
\newcommand{\Bv}{\bm{v}}
\newcommand{\Bd}{\bm{d}}
\newcommand{\BB}{\bm{B}}
\newcommand{\BW}{\bm{W}}
\newcommand{\Bxi}{\bm{\xi}}

\newcommand{\Btheta}{\bm{\theta}}
\newcommand{\tNN}{{\rm NN}}
\newcommand{\ta}{{\rm a}}
\newcommand{\tb}{{\rm b}}
\newcommand{\tn}{{\rm n}}

\newcommand{\dsigma}{\dot{\sigma}}
\newcommand{\hpNN}{\hat{p}_{{\rm NN}}}
\newcommand{\BRLU}{\mathcal{B}_\text{ReLU}}
\newcommand{\Bdsg}{\mathcal{B}_{\dot{\sigma}}}

\definecolor{orange}{RGB}{255,127,0}

\newcommand{\hz}[1]{\textcolor{blue}{#1}}

\newtheorem{lem}{Lemma}[section]
\newtheorem{theo}{Theorem}[section]

\newtheorem{assu}{Assumption}[section]

\title{Stationary Density Estimation of It\^o Diffusions Using Deep Learning}

\author{
Yiqi Gu \\
Department of Mathematics \\
National University of Singapore, 10 Lower Kent Ridge Road, Singapore, 119076 \\
({\texttt{matguy@nus.edu.sg})} \\
  \And
  John Harlim \\
  Department of Mathematics, Department of Meteorology and Atmospheric Science, \\ Institute for Computational and Data Sciences \\
  The Pennsylvania State University, University Park, PA 16802, USA\\
  \texttt{jharlim@psu.edu} \\
  \And
   Senwei Liang \\
  Department of Mathematics, Purdue University, IN 47907, USA \\
  \texttt{liang339@purdue.edu} \\
  \And
    Haizhao Yang \\
  Department of Mathematics, Purdue University, IN 47907, USA\\
  \texttt{haizhao@purdue.edu}
}

\begin{document}

\maketitle

\begin{abstract}
In this paper, we consider the density estimation problem associated with the stationary measure of ergodic It\^o diffusions from a discrete-time series that approximate the solutions of the stochastic differential equations. To take an advantage of the characterization of density function through the stationary solution of a parabolic-type Fokker-Planck PDE, we proceed as follows. First, we employ deep neural networks to approximate the drift and diffusion terms of the SDE by solving appropriate supervised learning tasks. Subsequently, we solve a steady-state Fokker-Plank equation associated with the estimated drift and diffusion coefficients with a neural-network-based least-squares method. We establish the convergence of the proposed scheme under appropriate mathematical assumptions, accounting for the generalization errors induced by regressing the drift and diffusion coefficients, and the PDE solvers. This theoretical study relies on a recent perturbation theory of Markov chain result that shows a linear dependence of the density estimation to the error in estimating the drift term, and generalization error results of nonparametric regression and of PDE regression solution obtained with neural-network models. The effectiveness of this method is reflected by numerical simulations of a two-dimensional Student's t distribution and a 20-dimensional Langevin dynamics.
\end{abstract}

\keywords{Stochastic differential equations \and Data-driven method \and Deep neural network \and Fokker-Plank equation}

\section{Introduction}\label{sec:intro}

Many phenomena subject to random perturbations can be modeled by stochastic differential equations~(SDEs) driven by Brownian noises. Under some regularity assumption, the time evolution of the probability measure can be characterized by the Fokker-Planck equation, a parabolic partial differential equation that depicts the time evolution of the density function of the underlying stochastic processes. Despite its wide applications in modeling physical or biological systems.~\cite{risken1996fokker,iancu2001nonlinear,frank2005nonlinear,chavanis2008nonlinear,hess1976fokker}, solving the Fokker-Planck PDE associated to high-dimensional It\^o diffusion processes is computationally a challenging task. In this paper, we are interested in estimating the density function associated with the stationary solution of the Fokker-Planck PDE from a discrete-time series of approximate solutions of the underlying SDEs without knowing the explicit drift and diffusion components.

Density estimation is a long-standing problem in computational statistics and machine learning.  Among the existing approaches, it is widely accepted that the classical Kernel Density Estimation (KDE) \cite{rosenblatt1956remarks} is not effective for problems with dimension higher than three (see e.g., \cite{hwang1994nonparametric, bsp:2013,wang2019nonparametric}). Along this line, the kernel embedding (another class of linear estimator) also suffered from the curse of dimension \cite{zhl_fods:2020}. Another class of popular parametric density estimators is the Gaussian Mixture Models (which is also known as the Radial Basis Models in some literature) \cite{hwang1994nonparametric}. This class of approaches is considered as a nonlinear estimator method since the training involves the minimization of a loss function that depends nonlinearly on the latent parameters. A practical issue of such a convex nonlinear optimization problem is the difficulty in identifying the global minimizer using numerical methods. While this issue is not solved, recent advances in deep learning theory show that the deep neural network (DNN), as a composition of multiple linear transformations and simple nonlinear activation functions, has the capacity of approximating various kinds of functions, overcoming or mitigating the curse of dimensionality~\cite{Hadrien,Weinan2019,bandlimit,poggio2017,yarotsky2019,Shen3,MO,HJKN19_814,Shen4}. Besides, it is shown that with over-parametrization and random initialization, the DNN-based least square optimization achieves a global minimizer by gradient descent with a linear convergence rate in both the setting of regression~\cite{DBLP:journals/corr/abs-1806-07572,Du2018,Zhu2019,Chen1,MeiE7665,lu2020meanfield,ding2021overparameterization,CHERIDITO2021101540} and PDE solvers~\cite{LuoYang2020,liang2021solving}. In parallel to this finding, several density estimators have adopted DNN, such as the Neural Autoregressive Distribution Estimation \cite{uria2016neural} and its variant, the Masked Autoregressive Flow \cite{papamakarios2017masked}.

Building on these encouraging results, we consider solving the density estimation problem where the target function is the density associated with the stationary measure of an It\^o process. With this prior knowledge, we propose to solve the density estimation problem following these two steps. First, we employ a deep learning algorithm to solve appropriate supervised learning tasks to uncover the drift and diffusion coefficients of the SDEs. Second, we solve the stationary Fokker-Planck PDE generated from the estimated drift and diffusion coefficients. While traditional grid-based numerical methods, such as finite element methods and finite difference methods~\cite{spencer1993numerical,kumar2006solution,sepehrian2015numerical} can be employed to solve the Fokker-Plank equation, they are usually limited to low-dimensional problems. On the other hand, neural network-based methods has been successfully used in solving high dimensional PDEs~\cite{Sirignano2018,Han2018,Liu2020Jul,Zang2020,RAISSI2019686,Khoo2017SolvingPP,zhai2020deep,Gu2020}, including the recent application in solving the high-dimensional Fokker-Plank equation~\cite{xu2020solving,zhai2020deep,liu2020neural}. These successes encourage us to also use deep learning to solve the approximate Fokker-Planck PDE.

We will also develop a new theory for the proposed approach with numerical verifications on low and relatively high-dimensional test examples, especially when the parameters of the Fokker-Planck equations have to be estimated, which has not been considered in the literature. Our theory can also explain and support the empirical success of existing deep learning approaches lacking the theoretical analysis of deep learning. The main goals of this theoretical study are to 1) understand under which mathematical assumptions can the density estimation problem be well-posed, 2) establish the convergence of the proposed scheme, and 3) identify the error in terms of training sample size, width/length of the neural-network models, discretization time step and noise amplitudes in the training data, and the dimension of the stochastic processes. In conjunction, we will also verify whether the perturbation theory \cite{zhang2021error} is valid. Particularly, we will check whether the stochastic process associated with the estimated drift and diffusion terms (obtained from deep learning regression in the first step) can indeed estimate the underlying invariant measure accurately. This verification is a  by-product that can practically be used to generate more samples if needed. 

The organization of this paper is as follows. In Section~\ref{sec:problemsetup}, we introduce the problem of stationary density estimation associated with It\^o diffusions. In Section~\ref{sec:1}, the deep learning method is discussed. In Section~\ref{sec:convergence}, we provide the convergence theoretical analysis. In Section~\ref{sec:numerical}, we present the numerical experiments of Student's distribution and Langevin dynamics. We conclude the paper with some remarks and open questions in Section~\ref{sec:cond}. To improve the readability, we report the proofs of the lemmas of Section~\ref{sec:convergence} in Appendix~\ref{appendixA}.

\section{Problem Setup}\label{sec:problemsetup}

Consider the following SDE,
\BEA
dX_t = \Ba(X_t)\,dt + \Bb(X_t)\,dW_t,
\label{eqn:sde}
\EEA
with an initial condition randomly drawn from an arbitrary well-defined distribution, $X_0 \sim \pi_0$. The SDE in \eqref{eqn:sde} is defined
with a drift term, $\Ba:\mathbb{R}^d\to\mathbb{R}^d$ and a diffusion tensor, $\Bb:\mathbb{R}^d\to \mathbb{R}^{d\times m}$, where $m\leq d$. Here, $W_t$ denotes the standard $m-$dimensional Wiener process. We assume that $\Ba$ and $\Bb$ are globally Lipschitz such that the SDE in \eqref{eqn:sde} with  the initial condition $X_0=x$ has a unique solution. In addition, we also assume that the Markov process $X_t$ is ergodic. This implies that the transition kernel corresponding to the Markov process $X_t$ converges to a unique stationary measure $\pi$ as $t\to\infty$. When the probability measure $\pi$ is absolutely continuous with respect to the Lebesque measure, $d\pi(x)= p(x)\,dx$, the density function $p:\mathbb{R}^d\to\mathbb{R}$ is the solution of the stationary Fokker-Planck equation,
\BEA
\OL^* p := -\mbox{div}(\Ba p) + \frac{1}{2} \sum_{i,j=1}^n \frac{\partial}{\partial x_i}\frac{\partial}{\partial x_j} ((\Bb\Bb^\top)_{ij} p) = 0,\label{eq:fpe}
\EEA
where $p\geq 0$ and $\int_{\mathbb{R}^d} p(x)\,dx = 1.$ We will state these (and additional) assumptions in Section~\ref{sec:convergence} for the convergence analysis study.

In this work, we aim to estimate the stationary density $p$ of the SDE \eqref{eqn:sde} without the knowledge of $\Ba$ and $\Bb$. What is available is a time series $\{\Bx^n\}_{n\geq 0}$ generated by a numerical SDE solver of \eqref{eqn:sde} that is assumed to possess an ergodic invariant measure, $\tilde{\pi}$, whose ``distance'' from $\pi$ can be controlled by the numerical discretization time step $\delta t$. We should point out that when $\Ba$ is globally Lipschitz and $\Bb$ is a full rank matrix and if the underlying Markov process in $X_t$ in \eqref{eqn:sde} is geometrically ergodic, then the Markov chain $\{\Bx^n\}$ induced by the Euler-Maruyama discretization is also geometrically ergodic \cite{mattingly2002ergodicity}. In Section~\ref{sec:convergence}, we will restrict our convergence study to this case.
In a less stringent case, e.g., $\Ba$ is locally Lipschitz,  the Markov chain induced by EM discretization is not ergodic in general. While one can generate an ergodic Markov chain by solving the SDE in \eqref{eqn:sde} with a stochastic backward Euler discretization \cite{mattingly2002ergodicity}, consistent learning from samples of such an ergodic chain will induce a more complicated loss function that incorporates the backward Euler scheme. While this case can be incorporated numerically, we neglect it in this paper since generally speaking the discretization scheme is unknown and the inconsistency of the numerical schemes that are used in generating the time series and in the construction of loss function in the learning algorithm induces an additional bias. For simplicity, we consider discrete Markov chain $\Bx^n$ generated by EM scheme,
\BEA
\Bx^{n+1} - \Bx^n = \Ba(\Bx^n)\delta t + \Bb(\Bx^n) \sqrt{\delta t} \Bxi_n, \quad\quad \Bxi_n\sim\mathcal{N}(0,\bm{I}_m),
\label{eqn:euler}
\EEA
where $\delta t$ denotes the time step size and $\bm{I}_m$ is an $m\times m$ identity matrix. In the next section, we will use the same discretization to construct the appropriate loss functions to approximate $\Ba$ and $\Bb\Bb^\top$. Since the available training data are sampled from $\tilde{\pi}$, the learning algorithm can only (at best) achieve a population risk defined with respect to $\tilde{\pi}$ and we will characterize the error induced by the EM discretization using an existing perturbation theory result.

While the SDE is defined on an entire unbounded domain $\mathbb{R}^d$ (the measure is not compactly supported or the density is strictly positive away from zero), numerically we can only solve the PDE on a bounded domain. Following existing approaches of solving Fokker-Planck PDEs with neural-networks \cite{xu2020solving,uy2020neural,zhai2020deep}, we consider a simply connected compact domain $\Omega\subset \mathbb{R}^d$ large enough such that the density on $ \mathbb{R}^d \backslash \Omega$ is effectively negligible. Practically, this assumption implies that the training data $\Bx^n \in \Omega$, and the stationary solution that we are looking for can be normalized with respect to $\Omega$, that is, $\int_\Omega p(\Bx)\,d\Bx = 1$. This assumption is critical especially when the vector field $\Ba$ is unknown and needs to be numerically estimated with deep learning, for which one can only (at best) guarantee the error in $L^2$- topology over a compact domain. In Section~\ref{sec:convergence}, we will clarify this assumption.


\section{Deep learning method for density estimation}\label{sec:1} In this section, we introduce a deep learning method to estimate the stationary density of SDE \eqref{eqn:sde} from a time series of its solution, which consists of two steps. We begin the discussion by reviewing two deep learning architectures that we will use in our numerical simulations, the fully connected neural network (FNN) and the residual neural network (ResNet) in Section ~\ref{sec:nn}. Given a time series of the SDEs in \eqref{eqn:sde}, we fit the drift $\Ba$ and diffusion coefficients $\Bb\Bb^\top$ in the SDE \eqref{eqn:sde} by NNs, denoted as $\Ba_\text{NN}$ and $\BB_\text{NN}$, respectively (see Section~\ref{sec:regression}). Define $\hat{\OL}^*$ as the Fokker-Planck (FP) differential operator generated from the estimated networks $\Ba_\text{NN}$ and $\BB_\text{NN}$ approximating the underlying (FP) operator $\OL^*$ in \eqref{eq:fpe}. Our approach in estimating the stationary density $p$ is to solve the homogeneous PDE $\hat{\OL}^*\hat{p}=0$, where $\hat{p}$ is a solution parameterized by an FNN. The PDE can be solved via the network-based least square method introduced in Section \ref{Sec_least_square}.

\subsection{Neural networks}\label{sec:nn}
We now give a brief overview of the two basic neural networks that have been widely employed in deep learning. The first one is the fully connected neural network (FNN). Suppose $d$ is the dimensions of inputs. Given an activation function $\sigma:\mathbb{R}\to\mathbb{R}$, $L\in\mathbb{N}^+$,  and $w_\ell\in\mathbb{N}^+$ for $\ell=1,\dots,L$, an FNN is constructed as the composition of $L$ simple nonlinear functions as follows
\begin{equation}
\phi_\text{NN}(\Bx;\Btheta):=\Bc^\top \Bh_{L} \circ \Bh_{L-1} \circ \cdots \circ \Bh_{1}(\Bx)\quad \text{for } \Bx\in\mathbb{R}^d,\notag
\end{equation}
where $\Bc\in \mathbb{R}^{w_L\times 1}$; $\bm{h}_{\ell}(\Bx_{\ell}):=\sigma\left(\BW_\ell \Bx_{\ell} + \Bg_\ell \right)$ with $\BW_\ell \in \mathbb{R}^{w_{\ell}\times w_{\ell-1}}$ and $\Bg_\ell \in \mathbb{R}^{w_\ell}$ for $\ell=1,\dots,L$ ($W_0:=d$). With the abuse of notations, $\sigma(\Bx)$ means that $\sigma$ is applied entry-wise to a vector $\Bx$ to obtain another vector of the same size. $w_\ell$ is the width of the $\ell$-th layer and $L$ is the depth of the FNN. $\Btheta:=\{\Bc,\,\BW_\ell,\,\Bg_\ell:1\leq \ell\leq L\}$ is the set of all parameters in $\phi_\text{NN}$ to determine the underlying neural network.

Besides FNN, in our numerical simulations, we will also consider the residual neural network (ResNet)~\cite{he2016deep}. Using similar notations above, ResNet can be defined recursively as follows,
\begin{align}
    &\bm{h}_0=x, \bm{h}_{-1}=\bm{0},\notag
    \\ &\mbv_\ell=\sigma\left(\BW_\ell \bm{h}_{\ell-1} + \Bg_\ell \right), \quad \ell=1,2,\cdots, L,\notag
    \\ &\bm{h}_\ell=\text{pad}(\bm{h}_{\ell-2})+\mbv_\ell, \quad \ell=1,2,\cdots, L,
    \label{skipconnection}
    \\ &\phi_\text{NN}(x;\Btheta)=\Bc^\top\bm{h}_L.\notag
\end{align}
Here, the function $\text{pad}(\cdot)$ is used to pad zeros to the vector such that two vectors in the summation~\eqref{skipconnection} are of same size.
Popular types of activation functions include the rectified linear unit (ReLU) $\sigma(x)=\max\{0,x\}$, ReLU$^3$ $\sigma(x)=\max\{0,x^3/6\}$, Tanh $\sigma(x)=\frac{e^{x}-e^{-x}}{e^{x}+e^{-x}}$ and Mish $\sigma(x)=x\text{Tanh}(\log(1+e^x))$~\cite{misra2019mish}.
We use $\mathcal{F}_{L,W,\sigma}$ to denote the class of FNNs with depth $L$, width $W$ for all layers and activation $\sigma$.

\subsection{Regression of drift and diffusion coefficients}\label{sec:regression}

Taking the expectation of ~\eqref{eqn:euler} with respect to $\Bxi_{n}$, one can see that
\BEA\label{01}
\mathbb{E} [\Bx^{n+1} - \Bx^n - \Ba(\Bx^n)\delta t] = 0.
\EEA
With this identity, we consider a supervised learning method for estimating $\Ba(\Bx)$ with neural networks. More precisely, we approximate every component of $\Ba(\Bx)$ by an FNN $a_\text{NN}(\Bx;\Btheta)$ parameterized by a set of trainable parameters $\Btheta$. In practice, letting $\By^n := \frac{\Bx^{n+1}-\Bx^n}{\delta t}$, by \eqref{01}, we define $\Btheta^{\ta}_i$ as follows,
\begin{equation}\label{eqn:traininga}
\Btheta^{\ta}_i := \arg\min_{\Btheta} \frac{1}{N}\sum_{n=0}^{N-1}\left|y_i^n - a_\text{NN}(\Bx^n;\Btheta) \right|^2,
\end{equation}
for $i=1,\cdots,d$, where $y_i^n$ is the $i$-th component of $\By^n$. Then we define the vector-valued function
\begin{equation}\label{eqn:estimator_a_NN}
\Ba_\text{NN}(\Bx;\Btheta^{\ta}):=\left[a_\text{NN}(\Bx;\Btheta^{\ta}_1),\cdots,a_\text{NN}(\Bx;\Btheta^{\ta}_d)\right]^\top
\end{equation}
as the drift estimator to approximate $\Ba(\Bx)$, where $\Btheta^{\ta}$ consists of $\{\Btheta^{\ta}_i\}$.

This is a supervised learning task to estimate $\Ba:\mathbb{R}^d\to \mathbb{R}^d$ from a pair of labelled training data set, $\{\Bx^n,\By^n\}_{n=0}^{N-1}$. To simplify the analysis in the next section, we assume that $\Bx^i$ are i.i.d. samples of the stationary random distribution $\tilde{\pi}$.
While we do not employ this simplification in our numerical study, practically, such i.i.d.~samples can be obtained by sub-sampling from the Markov chain $\{\Bx^n\}_{n\geq 0}$ such that their temporal correlation is negligible. For convenience of the following discussion, we denote $\mathcal{X}:=\{\Bx^0,\ldots, \Bx^{N-1}\}$ and $\mathcal{Y}:=\{\By^0,\ldots, \By^{N-1}\}$. In \eqref{eqn:traininga}, the parameter $\Btheta^{\ta}_i$ is a global minimizer of the empirical loss function. Practically, since tochastic gradient descent or the Adam method \cite{KingmaB14} is used, such a global minimizer may not necessarily be identified.

Next, we approximate $\Bb(\Bx)\Bb(\Bx)^\top$ in similar ways. The $(i,j)$-th component of $\Bb(\Bx)\Bb(\Bx)^\top$ can be approximated by an FNN $B_\tNN(\Bx;\Btheta^{\tb}_{ij})$. Since $\Bxi_n$ is independent of $\Bx^n$, using the fact $\mathbb{E}[\Bxi_n\Bxi_n^\top]=\bm{I}_n$ and \eqref{eqn:euler} we have
\BEA
\mathbb{E}\Big[(\Bx^{n+1} - \Bx^n - \Ba(\Bx^n)\delta t)(\Bx^{n+1} - \Bx^n - \Ba(\Bx^n)\delta t)^\top - \Bb(\Bx^n)\Bb(\Bx^{n})^\top\delta t  \Big] = 0.\notag
\EEA
Based on this identity, assuming that we have obtained the network $\Ba_\text{NN}(\Bx;\Btheta^{\ta})\approx\Ba(\Bx)$, we can compute $\Btheta^{\tb}_{ij}$ by
\begin{equation}\label{eqn:trainingbbt}
\Btheta^{\tb}_{ij}:=\arg\underset{\Btheta}{\min}\frac{1}{N}\sum_{i=0}^{N-1}\left|(\By^{n}_i-a_\text{NN}(\Bx^n,\Btheta^{\ta}_i))(\By^{n}_j-\Ba_\text{NN}(\Bx^n,\Btheta^{\ta}_j)))^\top-\frac{1}{\delta t}B_\text{NN}(\Bx^n;\Btheta) \right|^2.
\end{equation}
for $1\leq i,j\leq d$. Similarly, the global minimizer $\Btheta^{\tb}_{ij}$ may not be identified in practice. To summarize, If these global minimizers are identified, the training procedure gives $\BB_\text{NN}(\Bx):=\left[B_\text{NN}(\Bx,\Btheta^{\tb}_{ij})\right]_{i=1,\cdots,d}^{j=1,\cdots,d} \approx \Bb(\Bx)\Bb(\Bx)^\top$.

We should also point out that when the diffusion tensor is a constant matrix, $\Bb\in \mathbb{R}^{d\times m}$, we do not need to solve the optimization problem \eqref{eqn:trainingbbt} by deep learning. In such a case, $\BB_\text{NN}$ is specified as a matrix and we will empirically estimate $\Bb\Bb^\top$ using the residual from the drift estimator $\Ba_\text{NN}(\cdot)$. Particularly,
\BEA
\BB_\text{NN}:= \frac{\delta t}{N} \sum_{n=1}^{N} \left(\By^n - \Ba_\text{NN}(\Bx^n;\Btheta^{\ta}) \right) \left(\By^n - \Ba_\text{NN}(\Bx^n;\Btheta^{\ta}) \right)^\top,\label{constantBBtop}
\EEA
where we used the same notation $\BB_\text{NN}$ and understand that it is a $d\times d$ matrix in this case.

\subsection{Estimation of the stationary density}\label{Sec_least_square}
Given the approximate drift $\Ba_\text{NN} \approx \Ba$ and diffusion coefficients, $\BB_\text{NN} \approx \Bb\Bb^\top$, we define the estimated FP operator,
\BEA\label{eqn:estimatedFP}
\hat{\OL}^*p:= -\mbox{div}(\Ba_\text{NN}p) + \frac{1}{2} \sum_{i,j=1}^d \frac{\partial}{\partial x_i}\frac{\partial}{\partial x_j} (B_\text{NN}^{ij} p),
\EEA
where $B_\text{NN}^{ij}$ is the $(i,j)$-entry of $B_\text{NN}$.

Subsequently, the stationary density is estimated by solving the approximate stationary FP equation,
\BEA\label{02}
\hat{\OL}^*\hat{p}=0,\quad\text{in}~\Omega
\EEA
where $\hat{p}:\Omega \to (0,\infty)$ denotes the analytical solution of this PDE that satisfies,
\BEA
\int_{\Omega}\hat{p}(\Bx)d\Bx = 1.\label{03}
\EEA
Numerically, we set $\Omega$ to be a rectangular domain that is large enough yet tightly covers most of the data points in $\mathcal{X}$.

We solve the equation \eqref{02} with the condition \eqref{03} by the popular network-based least square method~\cite{doi:10.1002/cnm.1640100303,712178}. Specifically, We use a neural network $\hpNN(\Bx;\Btheta)$ with a parameter set $\Btheta$ determined by solving the following minimization problem,
\begin{equation}
\underset{\Btheta}{\min}~J[\hpNN(\cdot;\Btheta)],\notag
\end{equation}
where
\begin{equation}\label{eqn:loss}
J[q]:=\|\hat{\OL}^*q\|_{L^2(\Omega)}^2+\lambda_1 \left|\int_{\Omega} q(\Bx) d\bm{x} - 1\right|^2 + \lambda_2 \|q \|_{L^2(\partial \Omega)}^2,\quad\forall q:\Omega\rightarrow\mathbb{R},
\end{equation}

where $\lambda_1$ is a regularization constant corresponding to the normalization factor in \eqref{03} such to ensure nontrivial solution; $\lambda_2$ is a regularization parameter corresponding to an artificial Dirichlet boundary condition.
In our numerical simulation, we empirically found that the artificial boundary constraint can be neglected if the function values at the prescribed boundary is sufficiently small.

\comment{
In \eqref{eqn:loss}, we have defined,
\BEA\label{04}
\|f\|^2_{L^2(\Omega, \nu)} := \int_{\Omega} f(\Bx)^2 \td\nu(\Bx).
\EEA
}

In the practical computation, when $d$ is moderately large, the first term of \eqref{eqn:loss} is usually computed via a Monte-Carlo integration. For example, if the data $\{\Bx^n_{\rm I}\}_{n=1}^{N_1}$ are uniformly distributed points in $\Omega$, then
\BEA 
\|\hat{\OL}^*q\|_{L^2(\Omega)}^2 \approx \frac{|\Omega|}{N_1} \sum_{n=1}^{N_1} \Big|\hat{\OL}^*q(\Bx^n_{\rm I})\Big|^2,\label{eqn:diffeqn}
\EEA
where $|\Omega|$ denotes the volume of the domain $\Omega$.

Similarly, as for the second term in \eqref{eqn:loss}, Monte-Carlo integral is formulated as
\BEA\label{05}
\int_{\Omega} q(\Bx) d\Bx \approx \frac{|\Omega|}{N_2}\sum_{n=1}^{N_2} q(\Bx_{\rm II}^n),\label{eqn:moment2}
\EEA
where $\{\Bx_{\rm II}^n\}_{n=1}^{N_2}$ are uniformly distributed sampled points in $\Omega$.

For the third term in \eqref{eqn:loss}, we approximate,
\BEA
\|q\|_{L^2(\partial\Omega)} \approx \frac{|\partial\Omega|}{N_3}\sum_{n=1}^{N_3} |q(\Bx^n_{\rm III})|^2, \label{eqn:moment3}
\EEA
where $\{\Bx^n_{\rm III}\}_{n=1}^{N_3}$ are uniformly distributed sampled points in $\partial\Omega$.

Combining \eqref{eqn:diffeqn}, \eqref{eqn:moment2}, and \eqref{eqn:moment3}, the training procedure is to minimize the following empirical loss function,
\BEA
J_S[q] :=  \frac{|\Omega|}{N_1}\sum_{n=1}^{N_1} \left|\OL^* q(\Bx_{\rm I}^n)\right|^2 + \lambda_1 \left|\frac{|\Omega|}{N_2}\sum_{n=1}^{N_2} q(\Bx_{\rm II}^n)-1\right|^2 + \lambda_2 \frac{|\partial\Omega|}{N_3}\sum_{n=1}^{N_3} \left|q(\Bx_{\rm III}^n)\right|^2.\label{pdeempiricalloss}
\EEA
Let
\BEA\label{equ:trainingp}
\Btheta^S= \arg\min_{\Btheta}J_S[\hpNN(\cdot,\Btheta)],
\EEA
then the density estimator is given by $\hpNN(\cdot;\Btheta^S)\approx p(\cdot)$ with $\hpNN:\Omega \to \mathbb{R}$  and $\int_\Omega \hpNN(\Bx;\Btheta^S)\td \Bx \approx 1$.

We should point out that in our numerical simulations, since the time series $\{\Bx^n\}_{n=1}^N$ that are distributed in accordance to $\tilde{\pi}$ are available, we conveniently replace the first component in the loss function in  \eqref{eqn:loss} with a weighted norm, $L^2(\Omega,\tilde{\pi})$ and accordingly adjust the Monte-Carlo sum in the first component in the empirical loss function in \eqref{pdeempiricalloss}. While  the convergence analysis corresponding to a weighted norm is equivalent to that of the unweighted norm when $\tilde{pi}$ is absolutely continuous with respect to Lebesque measure with bounded density function, for simplicity of the exposition, we will consider the analysis corresponding to loss functions in \eqref{eqn:loss} with unweighted $L^2(\Omega)$ norms. If the dimension $d$ is lower, one can also adopt numerical quadrature rules such as Gauss-type quadrature to evaluate the integrals in \eqref{eqn:loss} for higher accuracy.

\section{Convergence Theory}\label{sec:convergence}
In this section, we deduce an error bound for the estimator $\hpNN(\Bx;\Btheta^S)$, where $\Btheta^S$ is the global minimizer of the empirical loss function in \eqref{pdeempiricalloss}. Throughout the discussion in this section, we restrict the diffusion coefficient $\Bb \in \mathbb{R}^{d\times m}$ to be a full column rank matrix. We use the notation $\|\cdot\|$ for the Euclidean norm in $\mathbb{R}^d$.

\subsection{Preliminary remarks}
Let us set the stage for our discussion by specifying the class of FNNs. In Section \ref{sec:nn}, we introduced the general class of FNNs $\mathcal{F}_{L,M,\sigma}$. While for the simplicity of analysis, we choose special classes of FNNs as the hypothesis spaces of the optimization.

On one hand, we consider using deep ReLU FNNs with uniform bounds in the minimization \eqref{eqn:traininga}, the regression of true drift $\Ba(\Bx)$. Specifically, for any $P>0$, we denote
\BEA\label{deep_FNN_bound}
\mathcal{F}_{L,M,\text{ReLU}}^P=\left\{ \phi\in\mathcal{F}_{L,M,\text{ReLU}}:~|\phi(\Bx)|\leq P,~\forall\Bx\in\Omega\right\},
\EEA
as the class of ReLU FNNs with depth $L$, width $M$, and a uniform bound $P$ in $\Omega$.

On the other hand, we consider using two-layer ReLU$^3$ FNNs with parameter bounds in the minimization \eqref{equ:trainingp}, the approximation of the true density $p(\Bx)$. More precisely, for any $Q>0$, we explicitly specify
\BEA\label{two_layer_FNN}
\mathcal{F}_{2,M,\dsigma,Q}=\left\{\phi:\Omega\rightarrow\mathbb{R}:~\phi(\Bx)=\frac{1}{M}\sum_{m=1}^Mc_m\dsigma(\Bw_m^\top\Bx),~ |c_m|,\|\Bw_m\|_1\leq Q\right\},
\EEA
where $\dsigma=\max(0,x^3/6)$ denoting the ReLU$^3$ activation function widely used in network-based methods for second-order PDEs. For simplicity, we omit the biases $\Bg_\ell$ in the definition of FNNs in Section \ref{sec:nn}.

Since the analysis depends on the results of the perturbation theory on the ergodic It\^o diffusion in \cite{zhang2021error}, we will briefly review the concepts of geometric ergodicity and other relevant results.

We will now make precise the assumptions mentioned in Section~\ref{sec:problemsetup}.

\begin{assu}\label{assumpdyn}
The following are key assumptions of the underlying system that generates the process $X_t$:
\begin{enumerate}
\item[i.] \textbf{Lipschitz \& Linear growth bound:} The vector field $\Ba:\mathbb{R}^d\to\mathbb{R}^d$ is globally Lipschitz with Lipschitz constant $\lambda_{\Ba}>0$ to ensure the existence and uniqueness of the solution of the SDE in \eqref{eqn:sde} given an initial condition.
There exists a constant $K\in(0,+\infty)$ such that
    \begin{equation*}
    \|\Ba(\Bx)\|^2\leq K^{2}(1 + \|\Bx\|^2), \forall \Bx\in\mathbb{R}^d.
    \end{equation*}
   This linear growth assumption will ensure that the even order moments can be bounded under the same rate.
   \comment{
   In our application below, consider $\hat{X}_n$ to be the discrete Markov chain induced by Euler-Maruyama scheme in \eqref{eqn:euler} at time step $\delta t$, then one can show that, for any $\ell\geq 1$,
   \BEA
   \mathbb{E}^{x}\Big[\|\hat{X}_n\|^{2\ell}\Big]\leq RD^n(1+\delta t)^n(1+\|x\|^{2\ell}),\label{momentdiscretebound}
   \EEA
for some constants $R, D >0$ that are independent of $\delta t>0$.}
   \item[ii.] \textbf{Geometric ergodicity:} The Markov process $X_t$ is geometrically ergodic with a unique invariant measure $\pi$. See e.g. Assumptions 2.2-2.3 in \cite{zhang2021error} for the detailed conditions to achieve the geometric ergodicity for the SDE driven by additive Brownian noises. One of the conditions that is important for our discussion is that there exists a Lyapunov function $V:\mathbb{R}^d\to [1,\infty)$ with $\lim\limits_{x\rightarrow \infty} V(\Bx)= +\infty$, and $c_1,c_2\in(0,+\infty)$ such that
\begin{equation*}
    \OL V(\Bx) \leq -c_1 V(\Bx) + c_2, \quad \forall \Bx\in \mathbb{R}^{d},
\end{equation*}
   where $\OL$ is the $L^2(\mathbb{R}^d)$ adjoint of the FP operator $\OL^*$ defined in \eqref{eq:fpe}.
 \item[iii.] \textbf{Essentially quadratic:} The Lyapunov function $V = W^{\ell}$ for some $\ell\geq 1$, where $W$ is essentially quadratic, i.e., there exist constants $C_{i}\in(0,+\infty)$, $i=1,2,3$, such that
\begin{equation} 
    C_{1}\left(1+\|\Bx\|^2 \right) \leq W(\Bx) \leq C_{2}\left(1 + \|\Bx\|^{2} \right), \quad \|\nabla W(\Bx)\|\leq C_{3}\left( 1+ \|\Bx\| \right), \quad \forall \Bx\in \mathbb{R}^{d}.\notag
\end{equation}
 Together with the previous two assumptions, there exists $\delta_0>0$ such that for all $\delta t \in (0,\delta_0)$, the discrete Markov chain induced by the EM algorithm in \eqref{eqn:euler} is geometrically ergodic with the invariant measure, $\tilde{\pi}$, and that,
\BEA
\sup_{f\in \mathcal{G}_{\ell}}\left|\pi(f)- \tilde{\pi}(f)\right| \leq K_1 (\delta t)^{\nu} \pi(V),\notag
\EEA
for some $K_1=K_1(\ell)$ and $\nu\in (0,1/2)$. Here, the supremum is defined over a set of locally Lipschitz functions bounded above by $V$,
\begin{equation}
\mathcal{G}_{\ell}: = \left\{f(\Bx)\leq V(\Bx), \forall \Bx\in \mathbb{R}^d \mbox{ and } \big| \;|f(\Bx)-f(\By)| \leq C_{\ell}\left(1+ \|\Bx\|^{2\ell-1} + \|\By\|^{2\ell-1}\right)\|\Bx-\By\|, \quad \forall \Bx,\By\in \mathbb{R}^{d} \right\}.\label{Gell}
\end{equation}
\end{enumerate}
\end{assu}

\begin{lem}\label{lemma_perterror} Under the assumptions~\ref{assumpdyn}, for any small $0<\epsilon \ll 1$, suppose that the estimator $\hat{\Ba}:\mathbb{R}^d\to\mathbb{R}^d$ is globally Lipschitz with Lipschitz constant independent of $\epsilon$ and is a consistent estimator in the following sense,
\BEA
\|\Ba(\Bx)- \hat\Ba(\Bx) \|^2 \leq K_2 (1+\|\Bx\|^2)\epsilon^2, \quad\forall \Bx\in \mathbb{R}^d,\label{consistentestimator}
\EEA
for some constant $K_2>0$ that is independent of $\epsilon$. Let us denote $\hat{X}_n :=\hat{X}(t_n)$, where $t_n = n\delta t$ to be a Markov chain generated by the solution to,
\BEA
d\hat{X} = \hat\Ba(\hat{X}_t)\,dt + \hat\Bb\,dW_t, \quad \hat{X}_0 = \Bx,\label{perturbedSDE}
\EEA
with $\hat\Bb\hat\Bb^\top:=\hat\BB$. For any $\Bx\in\mathbb{R}^d$, there exists $0<\rho<1$ and $K_1>0$ such that,
\BEA
\sup_{f\in \mathcal{G}_\ell}|\pi(f)-  \mathbb{E}^{\Bx}[f(\hat{X}_{n})] |  \leq K_3 \left[\left(\rho^{n} + \frac{1-\rho^{n}}{1-\rho}\epsilon\right)V(\Bx) \right], \quad \forall n\geq 0,\label{perturbbound}
\EEA
where the set $\mathcal{G}_\ell$ is defined in \eqref{Gell}. If the process $\hat{X}$ associated to \eqref{perturbedSDE} has an invariant measure $\hat{\pi}$, then there exist $0<\alpha<1$, $0<\beta<\infty$, and $0<\gamma<1-\alpha$ such that, $\hat{\pi}(V) \leq \frac{\beta}{1-\alpha-\gamma}$.
\end{lem}

The result above holds for all $\Bx \in \mathbb{R}^d$ by requiring the condition in \eqref{consistentestimator} and that underlying process $X(t)$ is ergodic in $\mathbb{R}^d$ with a unique invariant measure $\pi$. Similar conclusion was reported in \cite{huggins2017quantifying} under a much stronger uniform convergence in placed of \eqref{consistentestimator}. One of the key issue in applying this result directly to the learning configuration is that the assumption in \eqref{consistentestimator} can be difficult to achieve unless if one consider learning with a loss function defined with the topology that is used to deduced the error bound in \eqref{perturbbound}, which relies on the perturbation theory of Markov chain. The usual practical machine learning computations solve a supervised learning problem induced by a weaker topology (commonly $L^2$) on a bounded domain. In such a weaker topology (relative to the sup norm in \eqref{perturbbound}), one can at best expect to construct an estimator with convergence guaranteed under an $L^2(\Omega,\tilde{\pi})$ error on a compact domain $\Omega \supset \mathcal{X}$ that contains all the training data. In the numerical section, we will empirically show that the pointwise accuracy of $\Ba$ and verify the accuracy of the invariant mean and covariance statistics induced by a Markov chain generated by the estimated drift and diffusion coefficients.

To overcome the incompatibility of the domains, we consider the following assumption.

\begin{assu}\label{assumpdomain}
Let $\Omega\subset\mathbb{R}^d$ be a simply connected compact domain such that $P(X\notin\Omega)\leq \epsilon_0$ for some $0<\epsilon_0 \ll 1$. For example, let $\Omega:= B(0,R) = \{\Bx\in\mathbb{R}^d: \|\Bx\|\leq R\}$ be a closed Euclidean ball of radius $R>1$ and suppose that $X$ has mean zero (centered) and is a sub-exponentially distributed random variable, $SE(\nu^2,\alpha)$, with $\nu, \alpha>0$, then by concentration inequality for sub-exponential distribution, one obtains
\BEA
\mathbb{P}(\|X\|\geq R) \leq 2e^{-\frac{R}{2\alpha}} :=\epsilon_0, \quad\forall R>\nu^2\alpha^{-1}. \label{concentrationinequality}
\EEA
Let $\tilde{X}$ be a random variable corresponding to the stationary distribution induced by the Euler-Maruyama discretization in \eqref{eqn:euler}, using the Markov inequality and strong error bound of EM scheme, one can deduce that
$\mathbb{P}[\|X - \tilde{X}\| \leq (\delta t)^{1/4}]\leq (\delta t)^{-1/4} \mathbb{E}[|X - \tilde{X}] \leq C (\delta t)^{1/4}$, which means that $\mathbb{P}[\|\tilde{X}\|\geq R+(\delta t)^{1/4}] \leq \mathbb{P}(\|X\|\geq R) \mathbb{P}[\|X - \tilde{X}\| \geq (\delta t)^{1/4}] \leq O(\epsilon_0)$. Even if $X$ (resp. $\tilde{X}$) is defined on $\mathbb{R}^d$, one can almost surely realize $\|X\|\leq R$ (resp. $\|\tilde{X}\|\leq R+\delta t^{1/4}$) for large enough $R>0$. This assumption effectively means that the process $X$ satisfies the Assumption~\ref{assumpdyn} for $\Bx\in \Omega = B(0,R)$ almost surely for large enough $R$. This also implies that Lemma~4.1 is valid for $\Bx\in \Omega$, where we understood $\pi(f) := \int_{\Omega} f(\Bx) \pi(d\Bx)$ in \eqref{perturbbound}. In the convergence theory below, without loss of generality, we will assume that $\Omega=[0,1]^d$. For general $\Omega$, similar results can be derived easily by rescaling $\Omega$ to $[0,1]^d$ with an isomorphic map.
\end{assu}

With the above assumption, we only need to restrict our attention to a compact domain $\Omega$ and, hence, the assumption that $\hat{\Ba}$ is globally Lipschitz with Lipschitz constant independent of $\epsilon$ is reasonable. In our algorithm, we use the ReLU activation functions to construct  $\hat{\Ba}$ and, hence, $\hat{\Ba}$ is a globally Lipschitz continuous function. By  the simultaneous approximation of ReLU neural networks in \cite{GUHRING2021107,hon2021simultaneous}, as long as $\Ba\in C^s$ with $s>1$, there exists a ReLU network $\hat{\Ba}$ approximating $\Ba$ in the Sobolev norm of $W^{1,\infty}(\Omega)$ with $\Omega$ as a compact set. This means that the Lipschitz constant of $\hat{\Ba}$ can be bounded by a constant depending on $\Ba$ instead of the approximation accuracy. However, how to identify $\hat{\Ba}$ satisfying these assumptions is a problem of the optimization algorithm.

We use the notations $\tilde{\pi}(f) = \int_{\Omega} f(\Bx) d\tilde\pi(\Bx)$ and $\hat{\pi}(f) = \int_{\Omega} f(\Bx) d\hat\pi(\Bx)$ for integrals over $\Omega$.
With Assumption \ref{assumpdomain}, we now let the solution $\hat{p}:\Omega\to (0,\infty)$ of the approximate FP equation be the density of $\hat{\pi}$, defined with respect to the Lebesque measure, $d\hat{\pi} = \hat{p}(x) dx$. Since the PDE in \eqref{02} is defined with the estimated coefficients, namely, $\Ba_\text{NN}:\Omega\to\Omega$ as defined in \eqref{eqn:estimator_a_NN} and $\BB_\text{NN}\in \mathbb{R}^{d\times d}$ as defined in \eqref{constantBBtop}, the error analysis below will need to account for the errors induced by these estimations. Recall that $\Bb$ is a constant matrix and $\Ba_\text{NN}$ is the best empirical estimator from the chosen hypothesis space (e.g., a class of FNN-functions of the chosen architecture), obtained by regressing the labeled training data $\{\Bx^i,\By^i\}_{i=1}^N$, where $\Bx^i \in \mathcal{X}$ and $\By^i:= \Ba(\Bx^i)+ \bm{\eta}^i$, $\bm{\eta}^i \sim \mathcal{N}(\bm{0},(\delta t)^{-1}\Bb\Bb^\top)$.

To quantify the error of the diffusion estimator, one can subtract $\Bb\Bb^\top$ from the empirical estimator defined in \eqref{constantBBtop} and derive the following upper bound,
\BEA
\|\Bb\Bb^{\top} - \BB_\text{NN} \|_{2} \leq  \left\|\sum_{i=1}^{N} D_i   \right\|_2 + \delta t\, \mathbb{E}_{\tilde{\pi}}\left[\left\|(\Ba(X) - \Ba_\text{NN}(X;\Btheta^{\ta}))\right\|^2\right], \label{spectralerror}
\EEA
where for each $i=1,\ldots, N$,
\BEA
D_i :=\frac{\delta t}{N}\big(\By^i - \Ba_\text{NN}(\Bx^i;\Btheta^{\ta})\big)\big(\By^i - \Ba_\text{NN}(\Bx^i;\Btheta^{\ta})\big)^\top -\frac{1}{N} \Big(\delta t \mathbb{E}_{\tilde{\pi}}[(\Ba(X) - \Ba_\text{NN}(X;\Btheta^{\ta}))(\Ba(X) - \Ba_\text{NN}(X;\Btheta^{\ta}))^\top] + \Bb\Bb^\top\Big),\label{rvD}
\EEA
is an independent, random, symmetric matrix of mean zero. Since $\Bx^i$ is bounded almost surely, one can bound $D_i$ almost surely with large enough $R>0$. In such a case,  one can use a matrix concentration inequality to bound the first term in \eqref{rvD} with large enough training sample $N$.
Particularly, using the Matrix Bernstein inequality (e.g., Theorem~1.6.2 in \cite{tropp2012user}), if we define
\BEA
\epsilon:=\delta t\, \mathbb{E}_{\tilde{\pi}}\left[\left\|(\Ba(X) - \Ba_\text{NN}(X;\Btheta^{\ta})\right\|^2\right] \label{def_epsilon},
\EEA
and denote $\|D_i\|\leq \frac{D}{N}$ for some $D>0$, then the first term in \eqref{spectralerror} is smaller than $\epsilon$ with probability $1-2d\exp(-\frac{\epsilon^2/2}{O(N^{-2})+DN^{-1}\epsilon/3})>0$. This means, one can bound $\|\sum_{i=1}^ND_i\|\leq \epsilon$ with high probability by choosing $N \geq C\epsilon^{-1}\log{2d}$. We can therefore conclude that the spectral error in \eqref{spectralerror} is of order-$\epsilon$, which is the generalization error rate as defined in \eqref{def_epsilon}.

We should point out that the result in Lemma~\ref{lemma_perterror} does not assume the ergodicity of the Markov process $\hat{X}(t)$ generated by the SDE in \eqref{perturbedSDE}.
Suppose that $\hat{X}(t)$ is generated with $\hat{\Ba} = \Ba_\text{NN}$ and $\hat{\Bb}\hat{\Bb}^\top = \BB_\text{NN}$ has an invariant measure $\hat{\pi}$ on $\Omega$. Integrating \eqref{perturbbound} with respect to $\hat{\pi}$, we obtain,
\BEA
\Big|\pi(f) - \hat{\pi}(f) \Big| &=&   \Big|\pi(f) - \int_\Omega f(\Bx)\hat{\pi}(d\Bx)\Big| =  \Big| \pi(f) -\int_{\Omega} \mathbb{E}^{\Bx}[f(\hat{X}_{n})]\hat{\pi}(d\Bx)\Big| \leq K_3\hat{\pi}(V)\epsilon \label{measureerrors}
\EEA
as $n\to\infty$. To obtain \eqref{measureerrors}, we have used \eqref{perturbbound}. With this background, the error bound for $\hpNN(\Bx;\Btheta^S)$ can be deduced by accounting the regression error of $a$ and the error from the proposed PDE solver,
\BEA
\Big|\pi(f) - \int_{\Omega} f(\Bx) \hpNN(\Bx;\Btheta^S)\,d\Bx \Big| &\leq & \big|\pi(f) -\int_{\Omega} f(\Bx)\hat{p}(\Bx)d\Bx \big| +\left|\int_{\Omega} f(\Bx) \big(\hat{p}(\Bx)- \hpNN(\Bx;\Btheta^S)\big)\,d\Bx\right| \notag\\ &\leq&
K_3\hat{\pi}(V)\delta t \underbrace{\mathbb{E}_{\tilde{\pi}}\left[\left\|(\Ba(X) - \Ba_\text{NN}(X;\Btheta^{\ta})\right\|^2\right]}_{(I)} + \|f\|_{L^2(\Omega)} \underbrace{\|\hat{p}-\hpNN(\cdot;\Btheta^S) \|_{L^2(\Omega)}}_{(II)},
\label{maininequality}
\EEA
where we have used \eqref{measureerrors}, \eqref{def_epsilon}, and the Cauchy-Schwartz inequality. In the next two subsections, we will bound the terms (I) and (II) in \eqref{maininequality}.

\subsection{Regression error for the drift estimator}\label{sec:regression_error}
Now let us consider the error in the regression of the drift coefficients, namely, the minimization problem \eqref{eqn:traininga}. We will derive the $L^2$ error with respect to $\tilde{\pi}$ between the estimator $\Ba_\text{NN}(\Bx;\Btheta^\ta)$ and the true drift function $\Ba(\Bx)$. For this purpose, given a class $\mathcal{F}$ of functions: $\Omega\rightarrow\mathbb{R}$, we denote its pseudo dimension by $\text{Pdim}(\mathcal{F})$, which is the largest integer $m$ for which there is some $(\Bx_1,\cdots,\Bx_m,y_1,\cdots,y_m)\in\Omega^m\times\mathbb{R}^m$ such that for any $(b_1,\cdots,b_m)\in\{0,1\}^m$, there exists $f\in\mathcal{F}$ satisfying $f(\Bx_i)>y_i\Leftrightarrow b_i=1~\forall i$. The prediction error analysis of FNNs have been studied in several papers, e.g., \cite{10.1214/19-AOS1875,JMLR:v21:20-002,Chen2019NonparametricRO,LuoYang2020,Jiao2021,lu2021priori,mishra2020estimates,duan2021convergence}. In particular, we introduce the following lemma concerning the prediction error of the FNN-based least square regression, which is studied in \cite{Jiao2021}.

\begin{lem}[\cite{Jiao2021}, Theorem 4.2]\label{lem01}
Let $f_0:[0,1]^d\rightarrow\mathbb{R}$ be a H{\"o}lder continuous function, i.e., there exist $\lambda\geq0$ and $\alpha\in(0,1]$ such that $|f_0(x)-f_0(y)|\leq\lambda\|\Bx-\By\|^\alpha$ for all $\Bx,\By\in[0,1]^d$. Suppose $\|f_0\|_{L^\infty([0,1]^d)}\leq P$ for some $P\geq1$. Let $\nu$ be a probability measure that is absolutely continuous with respect to the Lebesgue measure and a random variable $\Bx\sim\nu$. Let $\eta$ be a random variable satisfying $\mathbb{E}[\eta]=0$ and $\text{Var}[\eta]=\sigma^2$. Let $\{\Bx^n\}_{n=1}^N$ be $N$ independent and identically distributed samples of $\Bx$, and $y^n=f_0(\Bx^n)+\eta$ is the response with noise $\eta$ for each $n$. For any $I_1,I_2\in\mathbb{N}^+$, let
\begin{equation}
\Btheta^{f_0}:= \arg\min_{\Btheta} \frac{1}{N}\sum_{n=1}^N\left|y^n - f_\text{NN}(\Bx^n;\Btheta) \right|^2,\notag
\end{equation}
where $f_\text{NN}\in\mathcal{F}_{L,W,\text{ReLU}}^P$ having depth $L=12I_2+14$ and width $W=\max\{4d\lfloor I_1^\frac{1}{d}\rfloor+3d,12I_1+8\}$ for all hidden layers. Then the prediction error is given by
\BEA\label{28}
\mathbb{E}_\nu\left[|f_\text{NN}(\cdot,\Btheta^{f_0})-f_0|^2\right]\leq C[P^2WL(d+WL)\log(Wd+W^2L)(\log N)^3N^{-1}
+\lambda^2d(I_1I_2)^{-4\alpha/d}],
\EEA
for $N\geq\text{Pdim}(\mathcal{F}_{L,W,\text{ReLU}}^P)$, where $C$ is a constant that does not depend on $d$, $N$, $L$, $W$, $\lambda$, $\alpha$, $I_1$, $I_2$, $P$.
\end{lem}

In Lemma \ref{lem01}, the exponent of the error bound in \eqref{28} can be improved to be  dimension-independent if we assume $f_0$ is in Barron-type spaces, which are first studied in \cite{barron1993} and further developed in \cite{Weinan2019,E2020,lu2021priori,siegel2021improved,siegel2021characterization,chen2021representation,Caragea2020}. Here we follow the Barron space with respect to two-layer ReLU networks proposed in \cite{E2020}. Suppose $f:\Omega\rightarrow\mathbb{R}$ is a function having the following form,
\BEA
f(\Bx)=\int_{\mathbb{R}\times\mathbb{R}^d} c\max(\Bw^\top\Bx,0)\rho(\td c,\td \Bw)=\mathbb{E}_\rho[c\max(\Bw^\top\Bx,0)],\quad\Bx\in\Omega\notag
\EEA
for some probability measure $\rho$ on $\mathbb{R}\times\mathbb{R}^d$, then its Barron norm is defined by
\BEA
\|f\|_{\BRLU}=\underset{\rho\in P_f}{\inf}(\mathbb{E}_\rho|c|\|\Bw\|_1),
\EEA
where $P_f:=\left\{\rho:f(\Bx)=\mathbb{E}_\rho[c\max(\Bw^\top\Bx,0)]\right\}$. And the ReLU Barron space is defined by $\BRLU=\{f\in C^0:\|f\|_{\BRLU}<\infty\}$. Now we have the following result.
\begin{lem}\label{lem02}
Let $f_0:[0,1]^d\rightarrow\mathbb{R}$ such that $\|f_0\|_{\BRLU}\leq P$ and $\|f_0\|_{L^\infty([0,1]^d)}\leq P$ for some $P\geq1$. For the least square regression proposed in Lemma \ref{lem01}, we let $f_\text{NN}\in\mathcal{F}_{2,W,\text{ReLU}}^P$ for some $W\in\mathbb{N}^+$, Then the prediction error is given by
\BEA\label{28}
\mathbb{E}_\nu\left[|f_\text{NN}(\cdot,\Btheta^{f_0})-f_0|^2\right]\leq C\left[P^2W(d+W)\log(Wd+W^2)(\log N)^3N^{-1}
+\|f_0\|_{\BRLU}^2dW^{-1}\right],
\EEA
for $N\geq\text{Pdim}(\mathcal{F}_{2,W,\text{ReLU}}^P)$, where $C$ is a constant that does not depend on $d$, $N$, $W$, $f_0$, $P$.
\end{lem}

\begin{proof}
See Appendix~\ref{appendixA}.
\end{proof}

In our case, we set in the hypothesis of Lemma \ref{lem01} that $L=O(I_2)$ and $W=O(I_1)$ are both large integers. Combining with Lemma \ref{lem02}, the error estimation for the minimization problem \eqref{eqn:traininga} can be directly obtained.

\begin{lem}\label{thm00}
In addition to the Assumption \ref{assumpdyn}, we let $\tilde{\pi}$ be absolutely continuous with respect to the Lebesgue measure. Denote $P_{\Ba}=\max\{\|\Ba\|_{L^\infty(\Omega)},1\}$.
\begin{enumerate}
\item Let $L$ and $W$ be integers large enough, then the estimator $\Ba_\text{NN}$ defined in \eqref{eqn:estimator_a_NN} with components $a_\text{NN}\in\mathcal{F}_{L,W,\text{ReLU}}^{P_{\Ba}}$ satisfies
\BEA\label{29}
\mathbb{E}_{\tilde{\pi}}\left[|\Ba_\text{NN}-\Ba|^2\right]\leq C_{\Ba} \left(d^2WLN^{-1}+d(WL)^2N^{-1}+d^2(WL)^{-4/d}\right),
\EEA
for $N\geq\text{Pdim}(\mathcal{F}_{L,W,\text{ReLU}}^{P_{\Ba}})$;
\item Suppose all components of $\Ba$ are in $\BRLU$ with Barron norms no greater than $P_{\Ba}$. Let $W\in\mathbb{N}^+$, then the estimator $\Ba_\text{NN}$ defined in \eqref{eqn:estimator_a_NN} with components $a_\text{NN}\in\mathcal{F}_{2,W,\text{ReLU}}^{P_{\Ba}}$ satisfies
\BEA\label{29}
\mathbb{E}_{\tilde{\pi}}\left[|\Ba_\text{NN}-\Ba|^2\right]\leq C_{\Ba} \left(d^2WN^{-1}+dW^2N^{-1}+d^2W^{-1}\right),
\EEA
for $N\geq\text{Pdim}(\mathcal{F}_{2,W,\text{ReLU}}^{P_{\Ba}})$,
\end{enumerate}
where $C_{\Ba}>0$ is a term that depends on $\Ba$ and at most a polynomial in the logarithm of $N$, $L$, $W$.
\end{lem}

In Lemma \ref{thm00}, the Barron assumption on the target function helps to overcome the curse of dimensionality. In the following analysis for the solution error in the approximate FP equation, we will specify a Barron space for ReLU$^3$ networks and assume that the true solution is in this space; therefore the derived solution error is also exponentially independent of dimensions.

Another situation to mitigate the curse of dimensionality is when the data points are supported on a neighborhood of a low-dimensional Riemannian submanifold in $\Omega$ \cite{Chen2019NonparametricRO,Jiao2021}. Since it does not apply to the current problem in practice, we will not discuss more on this situation.

\subsection{Solution error for the approximate FP equation}
Now let us consider the error between $\hpNN(\cdot;\Btheta^S)$ and the true solution $\hat{p}$ of the approximate stationary FP equation \eqref{02}. In this section, we only consider the case that  $\{\Bx^n_{\rm MC}\}_{n=1}^{N_1}$ in \eqref{05} are uniformly distributed in $\Omega$. Similar results apply to other measures with smooth densities supported on $\Omega$.

First, we rewrite the approximate stationary FP equation (15) in the following divergence form
\begin{equation}\label{eqn:PDE_divergence_form}
-\hat{\OL}^*\hat{p}=-\sum_{i,j=1}^d\left(\frac{1}{2}B_\tNN^{ij}\hat{p}_{x_j}\right)_{x_i}+\sum_{i=1}^da_\tNN^i\hat{p}_{x_i}
+\left(\sum_{i=1}^d\frac{\partial a_\tNN^i}{\partial x_i}\right)\hat{p}=0,\quad\text{in}~\Omega.
\end{equation}

The error analysis is valid only when the equation \eqref{eqn:PDE_divergence_form} is well-posed. So we need to set up specific assumptions on the coefficients of \eqref{eqn:PDE_divergence_form}. First, note that $\BB_\tNN$ is positive semi-definite and \eqref{eqn:PDE_divergence_form} is elliptic, so we assume further that \eqref{eqn:PDE_divergence_form} is non-degenerate by specifying the smallest eigenvalue of $\BB_\tNN$ as a positive number. Also, we assume that the coefficients have a uniform bound, which is common in the analysis of elliptic equations.
\begin{assu}\label{assum_01}
The smallest eigenvalue of the symmetric matrix $\BB_\tNN$, denoted as $\Lambda$, is positive. Besides, $|B_\tNN^{ij}|<2B_1$, $|a_\tNN^i(\Bx)|<B_1$, $\left|\sum_{i=1}^d\partial a_\tNN^i(\Bx)/\partial x_i\right|<B_1$, $\forall i,j$ and $\forall\Bx\in\Omega$, for some $B_1>0$ .
\end{assu}

Next, considering \eqref{eqn:PDE_divergence_form} is defined in a compact domain, we can not guarantee the uniqueness of the solution $\hat{p}$ since no boundary condition is specified. Moreover, even if we impose a boundary condition, say Dirichlet condition $\hat{p}=g$ on $\partial\Omega$, we still need extra assumptions on the coefficients to ensure the uniqueness. For the latter, it suffices to take the following assumption.

\begin{assu}\label{assum_02}
\begin{equation*}
\int_\Omega\sum_{i=1}^da_\tNN^iv_{x_i}\cdot v+\left(\sum_{i=1}^d\frac{\partial a_\tNN^i}{\partial x_i}\right)v^2\td \Bx\geq 0,\quad\forall v\in H^1(\Omega).
\end{equation*}
\end{assu}

Under Assumption \ref{assum_02}, one can show that \eqref{eqn:PDE_divergence_form} with any Dirichlet condition admits a unique solution by Fredholm alternative and Lax-Milgram theorem. However, we can not specify such a boundary condition since no information on $\partial\Omega$ is provided. Fortunately, we note that the true density $p$ vanishes as $|x|\rightarrow\infty$, so it can be assumed that the approximate density $\hat{p}$ has a similar behavior. Although we do not specify any boundary value for $\hat{p}$, we can assume that $\hat{p}$ ``almost" vanishes on  $\partial\Omega$ as follows.

\begin{assu}\label{assum_03}
Let
$\|\hat{p}\|_{L^\infty(\partial\Omega)}\leq\epsilon_{\hat{p}}$ and $\|\hat{p}\|_{H^1(\partial\Omega)}\leq\epsilon_{\hat{p}}$ for some small positive number $\epsilon_{\hat{p}}>0$.
\end{assu}

Under Assumption \ref{assum_02} and \ref{assum_03}, it can be shown that any two solutions of \eqref{eqn:PDE_divergence_form} are close to each other up to accuracy $\epsilon_{\hat{p}}$ by standard elliptic equation analysis.

Now we indicate that the error $\|q-\hat{p}\|_{L^2(\Omega)}$ for any function $q$ is bounded by the loss function $J[q]$ and $\epsilon_{\hat{p}}$.

\begin{lem}\label{thm02}
Assume $\hat{p}$ is a classical solution of \eqref{02} with the condition \eqref{03}.
Let $q\in C^2(\bar{\Omega})$ and assume $\|\nabla q\|_{L^2(\partial\Omega)}\leq B_2$ for some $B_2>0$. If Assumptions \ref{assum_01}-\ref{assum_03} hold, then
\begin{equation*}
\|q-\hat{p}\|_{L^2(\Omega)}^2\leq C\left(J[q]+d(1+\epsilon_{\hat{p}})J[q]^\frac{1}{2}+d(1+\epsilon_{\hat{p}})\epsilon_{\hat{p}}\right),
\end{equation*}
where $C$ only depends on $\Omega$, $\Lambda$, $B_1$, $B_2$, $\lambda_1$, $\lambda_2$.
\end{lem}

\begin{proof}
See Appendix~\ref{appendixA}.
\end{proof}

Next, we estimate $J[\hat{p}]$ via the generalization analysis of FNNs. In the analysis, we redefine the Barron space for two-layer ReLU$^3$ networks and assume $\hat{p}$ is in this Barron space. The definition directly follows the ReLU Barron space proposed in Section \ref{sec:regression_error} except that we replace the ReLU activation with the ReLU$^3$ activation. Accordingly, we slightly modify the  Barron norm, which is also proposed in \cite{LuoYang2020}. Recall that $\dsigma$ denotes the ReLU$^3$ activation function, i.e. $\dsigma=\max(0,x^3/6)$.

Suppose $f:\Omega\rightarrow\mathbb{R}$ is a function having the following form,
\BEA
f(\Bx)=\int_{\mathbb{R}\times\mathbb{R}^d} c\dsigma(\Bw^\top\Bx)\rho(\td c,\td \Bw)=\mathbb{E}_\rho[c\dsigma(\Bw^\top\Bx)],\quad\Bx\in\Omega\notag
\EEA
for some probability measure $\rho$ on $\mathbb{R}\times\mathbb{R}^d$, then its ReLU$^3$ Barron norm is defined by
\BEA\label{07}
\|f\|_{\Bdsg}=\underset{\rho\in P_f}{\inf}(\mathbb{E}_\rho|c|\|\Bw\|_1^3),
\EEA
where $P_f:=\{\rho:f(\Bx)=\mathbb{E}_\rho[c\dsigma(\Bw^\top\Bx)]\}$. And the ReLU$^3$ Barron space is defined by $\Bdsg=\{f\in C^0:\|f\|_{\Bdsg}<\infty\}$. Now let us derive the uniform approximation of FNNs in $\mathcal{F}_{2,M,\dsigma,Q}$ for Barron functions.

\begin{lem}\label{lem03}
Given $f\in\Bdsg$, there exists some $p_\tNN\in\mathcal{F}_{2,M,\dsigma,\max\left\{\|f\|_{\Bdsg}/M,1\right\}}$ such that
\BEA\label{10}
\underset{\Bx\in\Omega}{\sup}~\left|\hat{\OL}^*p_\tNN(\Bx)-\hat{\OL}^*f(\Bx)\right|
+\underset{\Bx\in\Omega}{\sup}~\left|p_\tNN(\Bx)-f(\Bx)\right|
+\underset{\Bx\in\partial\Omega}{\sup}~\left|p_\tNN(\Bx)-f(\Bx)\right|\leq \left(4B_1+2\right)\|f\|_{\Bdsg}\sqrt{d/M},
\EEA
\end{lem}

\begin{proof}
See Appendix~\ref{appendixA}.
\end{proof}

Next, we introduce the error estimate for the Monte-Carlo integration, which can be directly proved using the Hoeffding's inequality.
\begin{lem}\label{lem04}
Given a compact domain $\Omega$. Suppose $f:\Omega\rightarrow\mathbb{R}$ is a function with $\|f\|_\infty<\infty$. Let $\{\Bx_n\}_{n=1}^N$ be a set of uniformly distributed points in $\Omega$. Then for any $\delta\in(0,1)$, with probability at least $1-\delta$ over the choice of $\Bx_n$,
\BEA
\left|\frac{|\Omega|}{N}\sum_{n=1}^Nf(\Bx_n)-\int_\Omega f(\Bx)\td\Bx\right|\leq\sqrt{\frac{2\|f\|_\infty^2\log(2/\delta)}{N}}.\notag
\EEA
\end{lem}

Now, the error estimate for the approximate FP equation is given as follows.

\begin{lem}\label{thm03}
Under Assumption \ref{assum_01}-\ref{assum_03}, further assume $\hat{p}\in\Bdsg$. Let $\Btheta^S=\text{argmin}_{\Btheta} J_S[\hpNN(\cdot,\Btheta)]$ with $\hpNN\in\mathcal{F}_{2,M,\dsigma,Q}$. Also, suppose $\{\Bx_\text{I}^n\}_{n=1}^{N_1}\subset\Omega$, $\{\Bx_\text{II}^n\}_{n=1}^{N_2}\subset\Omega$, $\{\Bx_\text{III}^n\}_{n=1}^{N_3}\subset\partial\Omega$ in \eqref{pdeempiricalloss} are uniformly distributed. Then for any $\delta\in(0,1)$, with probability of at least $1-\delta$ over the choice of these points,
\begin{equation}\label{25}
\|\hpNN(\Bx;\Btheta^S)-\hat{p}\|_{L^2(\Omega)}^2\leq C\left(J[\hpNN(\Bx;\Btheta^S)]+d(MQ^4d^\frac{1}{2}+\epsilon_{\hat{p}})J[\hpNN(\Bx;\Btheta^S)]^\frac{1}{2}+d(MQ^4d^\frac{1}{2}+\epsilon_{\hat{p}})\epsilon_{\hat{p}}\right),
\end{equation}
and
\begin{equation}
J[\hpNN(\Bx;\Btheta^S)]\leq C\left[I_1(Q,d,\delta,M,N_1,N_2,N_3)+I_2(Q,\delta,M,N_2)+I_3(\hat{p},d,\delta,M,N_2)\right],\notag
\end{equation}
with
\begin{align}
I_1&=(Q^8+1)\left(d^2\sqrt{\log(d)}+\log(Q^4+1)+\sqrt{\log(1/\delta)}\right)M^2(1/\sqrt{N_1}+1/\sqrt{N_3}),\notag\\
I_2&=MQ^4\sqrt{\log(6/\delta)/N_2}\left(MQ^4(\sqrt{\log(6/\delta)/N_2}+1)+1\right),\notag\\
I_3&=\|\hat{p}\|_{\Bdsg}^2d/M+\|\hat{p}\|_\infty^2\log(6/\delta)/N_2+\epsilon_{\hat{p}}^2,\notag
\end{align}
where $C$ only depends on $\Omega$, $\Lambda,B_1$, $\lambda_1$, and  $\lambda_2$. Especially, suppose $J[\hpNN(\Bx;\Btheta^S)]\leq1$ and let $N_\text{p}:=\min\{N_1,N_2,N_3\}$. Take $Q\leq O(M^{-\frac{1}{4}}d^{-\frac{1}{8}})$ and $N_\text{p}\geq O(\log(1/\delta))$ , then
\BEA
\|\hpNN(\Bx;\Btheta^S)-\hat{p}\|_{L^2(\Omega)}^2\leq O\left(d^2(\log(d))^\frac{1}{4}MN_\text{p}^{-\frac{1}{4}}+d^\frac{3}{2}M^{-\frac{1}{2}}+dN_\text{p}^{-\frac{1}{2}}+d\epsilon_{\hat{p}}\right),\notag
\EEA
with an order constant depending on $\Omega$, $\Lambda,B_1$, $\lambda_1$, $\lambda_2$, $\delta$, and $\hat{p}$.
\end{lem}

\begin{proof}
See Appendix~\ref{appendixA}.
\end{proof}

In Lemma \ref{thm03}, it implies using $N_\text{p}\sim O(M^s)$ with $s>4$ will reduce the solution error up to $O(d\epsilon_{\hat{p}})$ as $M,N_\text{p}\rightarrow\infty$. In practice, as an approximation of the original density $p$ which vanishes as $\|\Bx\|\rightarrow\infty$, the solution $\hat{p}$ could have a similar behavior. Hence $\epsilon_{\hat{p}}$ is small enough if $\Omega$ is moderately large. And this also leads to a small solution error $\|\hpNN-\hat{p}\|_L^2(\Omega)$.

Recall in the analysis of regression error for the drift estimator $\Ba$, we derive an error estimate for deep networks of any width and depth. While in the error analysis for the FP solution $\hat{p}$, only results for two-layer shallow networks are derived in the current work. It is promising to develop this analysis for deep networks in the future work.

\subsection{The main error estimation}
Inserting the two error bounds in Lemma~\ref{thm00} and Lemma~\ref{thm03} into the inequality in \eqref{maininequality} and collecting all the assumptions, we can show the following main theorem for the error estimation of the proposed algorithm.

\begin{theo}\label{thm01}
Let $\pi$ be the invariant measure of a Markov process $X_t$ that satisfies Assumptions~\ref{assumpdyn} and \ref{assumpdomain}. Let $P_{\Ba}\geq1$ such that $\|\Ba\|_{L^\infty(\Omega)}\leq P_{\Ba}$. Given discrete samples $\{\Bx^n\}_{n=0}^N$ of an ergodic measure $\tilde{\pi}$ that is absolutely continuous with respect to the Lebesque measure in $\mathbb{R}^d$, suppose that $\Ba_{NN}$ defined by \eqref{eqn:estimator_a_NN} with components $a_\text{NN}\in\mathcal{F}_{L,W,\text{ReLU}}^{P_{\Ba}}$ is a consistent estimator in the sense of \eqref{consistentestimator} for all $\Bx\in \Omega=[0,1]^d$. Suppose also that $N\geq\text{Pdim}(\mathcal{F}_{L,W,\text{ReLU}})$. Let the assumptions in Lemma~\ref{thm03} be valid, namely the Assumptions \ref{assum_01}-\ref{assum_03}. Suppose that $\hat{p}\in \Bdsg$ is estimated by $\hpNN(\cdot,\Btheta^S)\in\mathcal{F}_{2,M,\dsigma,Q}$ with $Q\leq O(M^{-\frac{1}{4}}d^{-\frac{1}{8}})$, where $\Btheta_S$ is the global minimizer of the empirical loss function \eqref{pdeempiricalloss}. Then, for all $f\in \mathcal{G}_\ell$ as defined in \eqref{Gell}, and for any $\delta\in (0,1)$, with probability of at least $1-\delta$ over the choice of $\{\Bx_\text{I}^n\}_{n=1}^{N_1}$, $\{\Bx_\text{II}^n\}_{n=1}^{N_2}$ and $\{\Bx_\text{III}^n\}_{n=1}^{N_3}$,
\BEA
\sup_{f\in\mathcal{G}_\ell}\Big|\pi(f) - \int_{\Omega} f(\Bx) \hpNN(\Bx,\Btheta^S)\,d\Bx \Big| &\leq& K_3\hat{\pi}(V)\delta t C_{\Ba} \left(d^2WLN^{-1}+d(WL)^2N^{-1}+d^2(WL)^{-4/d}\right) \notag \\ && + C_{\hat{p}}\left(d^2(\log(d))^\frac{1}{4}MN_\text{p}^{-\frac{1}{4}}+d^\frac{3}{2}M^{-\frac{1}{2}}+dN_\text{p}^{-\frac{1}{2}}+d\epsilon_{\hat{p}}\right), \label{mainerrorbound}
\EEA
where, $N_\text{p}:=\min\{N_1,N_2,N_3\}$ that satisfies $N_\text{p}\geq O(\log(1/\delta))$. Here, the term  $C_{\Ba}>0$ depends on $\Ba$ and at most a polynomial in the logarithm of $N$, $L$, $W$,
and the constant $C_{\hat{p}}>0$ depends on $\Omega$, $\delta$, $\hat{p}$, $\|f\|_{L^2(\Omega)}$, the regularization weights $\lambda_1$, $\lambda_2$, and the upper bounds constants $\Lambda, B_1$ defined in Assumption~\ref{assum_01}.
\end{theo}

And the error bound independent of dimensions exponentially is given as follows.
\begin{theo}
Under the hypothesis of Theorem \ref{thm01}, we further assume all components of $\Ba$ are in $\BRLU$ with Barron norms no greater than $P_{\Ba}$, and let $a_\text{NN}\in\mathcal{F}_{2,W,\text{ReLU}}^{P_{\Ba}}$, then the error bound term $d^2WLN^{-1}+d(WL)^2N^{-1}+d^2(WL)^{-4/d}$ in \eqref{mainerrorbound} can be improved to be $d^2WN^{-1}+dW^2N^{-1}+d^2W^{-1}$.
\end{theo}

We should point out that while the results are valid for the global minimizers $\Btheta_i^\ta$ in \eqref{eqn:traininga} and $\Btheta_S$ in \eqref{equ:trainingp}, we do not specify the condition for which such global minimizers are attainable. We directly assume that the minimizers are found, and do not consider the error from the optimization algorithms. In practice, one can not ensure that the global minimizers can be necessarily found by usual optimizers like gradient descent.


Moreover, throughout the convergence analysis, we consider using special FNN class \eqref{deep_FNN_bound} with uniform bounds or \eqref{two_layer_FNN} with parameter bounds as the hypothesis space, and derive corresponding approximation errors. However, in practical deep learning, one usually uses the general FNN class $\mathcal{F}_{L,W,\sigma}$ since it is closed under gradient descent optimizers and therefore easy for implementation.

\section{Numerical Examples}\label{sec:numerical}
In this section, we numerically demonstrate the effectiveness of our proposed methods on two test problems. The first example is a two-dimensional SDE with Student's t stationary distribution. The second example is a 20-dimensional Langevin dynamics associated to Lenard-Jones potential with Gibbs invariant measure.

In our examples, we replace the norm in the first term in \eqref{eqn:loss} with a weighted $L^2(\Omega,\tilde{\pi})$, recalling that $\tilde{\pi}$ denotes the stationary measure of the discrete Markov chain induced by \eqref{eqn:euler}. Hence we can directly use the available data set $\mathcal{X}:=\{\Bx^0,\ldots,\Bx^{N-1}\}$ as the Monte Carlo integration points. Empirically, we approximate the first term of \eqref{eqn:loss} via the following Monte-Carlo average,
\BEA
\|\hat{\OL}q\|_{L^2(\Omega,\tilde{\pi})}^2 \approx \frac{1}{|\mathcal{X}\cap\Omega|} \sum_{n=0}^{N-1} \Big|\hat{\OL} q(\bm{x}^n)\Big|^2 \mathbb{1}_{\Omega}(\Bx^n),
\EEA
where $\mathbb{1}_{\Omega}$ denotes the characteristic function over domain $\Omega$.

\subsection{Student's t-distribution}
Consider a two-dimensional SDE~\eqref{eqn:sde} for Student's t-distribution~\cite{averina1988numerical} with
\begin{align*}
    \Ba(\Bx)=\begin{bmatrix}
-\frac{3}{2}x_1+x_2 \\
\frac{1}{4}x_1-\frac{3}{2}x_2
\end{bmatrix}, \quad
\Bb(\Bx)= \begin{bmatrix}
\sqrt{\phi(x_1,x_2)} & 0 \\
-\frac{11}{8}\sqrt{\phi(x_1,x_2)}&\frac{\sqrt{255}}{8}\sqrt{\phi(x_1,x_2)}
\end{bmatrix}
\end{align*}
where $\Bx=(x_1,x_2)$ and $\phi(x_1,x_2) = 1+\frac{2}{15}(4x_1^2-x_1x_2+x_2^2)$. The stationary density is explicitly given by
\begin{align}
p(x_1,x_2) = \frac{2}{\pi\sqrt{15}}\left(\phi(x_1,x_2)\right)^{-3}.
\label{student solution}
\end{align}

\subsubsection{Data generation and implementation details} The time series dataset $\{\Bx^{i}\}_{i=0}^N$ is generated by EM scheme~\eqref{eqn:euler} with $\delta t= 0.05$ and $N=2\times10^7$. The bounded domain $\Omega$ is set as $[-4,4]\times[-6,6]$ such that over $98\%$ points are in $\Omega$.

 In our implementation, we use 6-hidden-layer ResNets (discussed in Section~\ref{sec:nn}) with the same width 50 per hidden layer and the smooth Mish activation function~\cite{misra2019mish}. To learn $\Ba_\text{NN}$ and $\Bb_\text{NN}$, Adam algorithm is applied to optimize the loss~\eqref{eqn:traininga} and~\eqref{eqn:trainingbbt} with batch size 10,000 for $T=20,000$ iterations. We use an initial learning rate of $10^{-4}$. The learning rate follows cosine decay with the increasing training iterations, i.e., the learning rate decays by multiplying a factor $0.5( cos(\frac{\pi t}{T}) + 1)$, where $t$ is the current iteration. To solve the PDE~\eqref{eqn:estimatedFP}, we optimize the loss~\eqref{eqn:loss} with $\lambda=1, \gamma=500$. In Adam, we use the batch size 10,000 for first term in~\eqref{eqn:loss} and 4,000 for the boundary term while the second term is approximated by $300^2$ Gaussian quadrature points. The learning rate is initialized by $10^{-3}$ and follows cosine decay.

\subsubsection{Identification of the drift and diffusion coefficients.}
To evaluate the accuracy, we define relative $L_2$ error as follows,
\begin{equation}
\frac{\|f-\hat{f}\|_{L^2(\Omega)}}{\|f\|_{L^2(\Omega)}},\label{relativeL2error}
\end{equation}
where $f$ and $\hat{f}$ represents the true function and the approximate function, respectively. Numerically, we approximate the integral over $10000$ Gaussian quadrature points in $\Omega$.

\begin{figure}
\centering
\subfloat[$\bm{a}_1$]{
\includegraphics[scale=0.30]{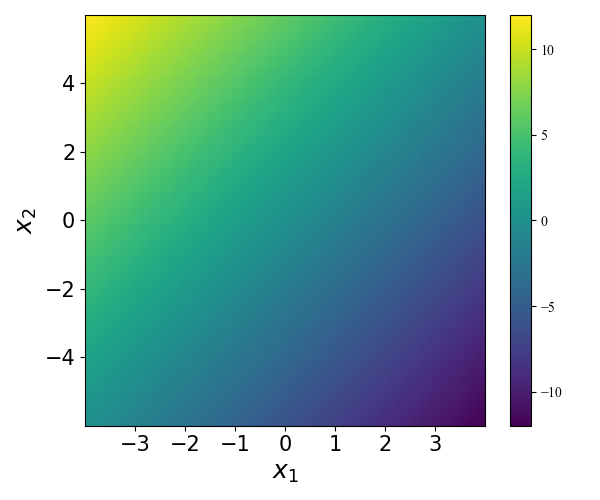}}
\hspace{0.1cm}
\subfloat[$(\Ba_\text{NN})_1$]{
\includegraphics[scale=0.30]{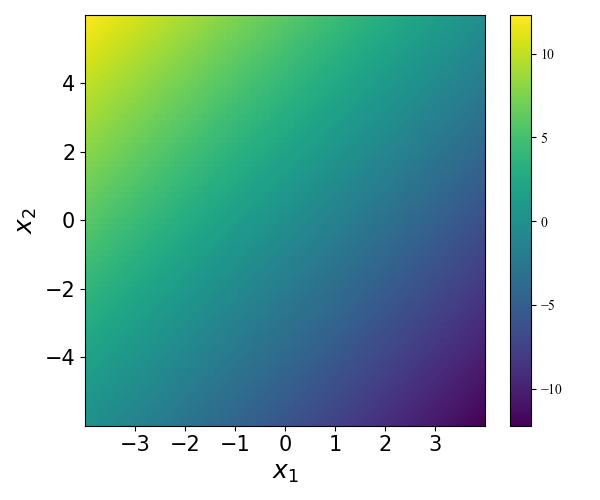}}
\hspace{0.1cm}
\subfloat[$\bm{a}_1$-$(\Ba_\text{NN})_1$]{
\includegraphics[scale=0.30]{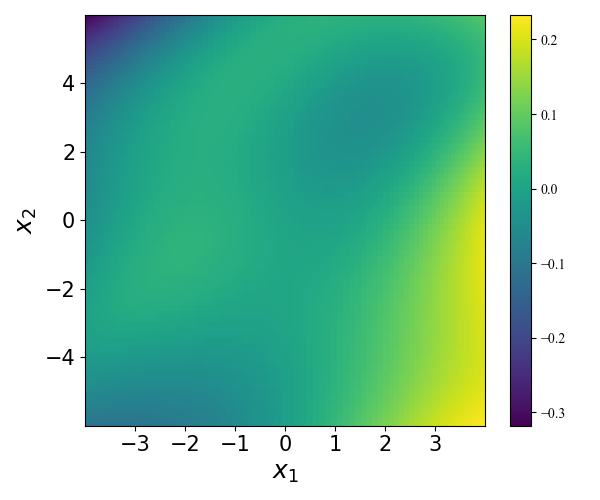}}
\caption{The comparison of the first component of drift term. (a) $\bm{a}_{1}$, (b) $(\Ba_\text{NN})_1$, and (c) their difference.}
\label{fig:a_approximate}
\end{figure}

\begin{figure}
\centering
\subfloat[$(\bm{b}\bm{b}^\top)_{11}$]{
\includegraphics[scale=0.30]{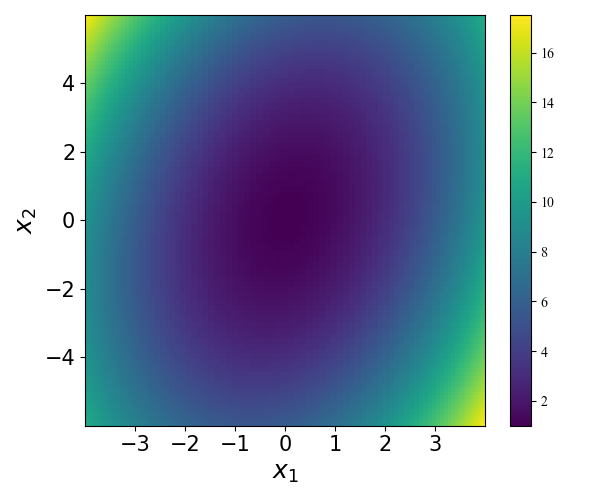}}
\hspace{0.1cm}
\subfloat[$(\BB_\text{NN})_{11}$]{
\includegraphics[scale=0.30]{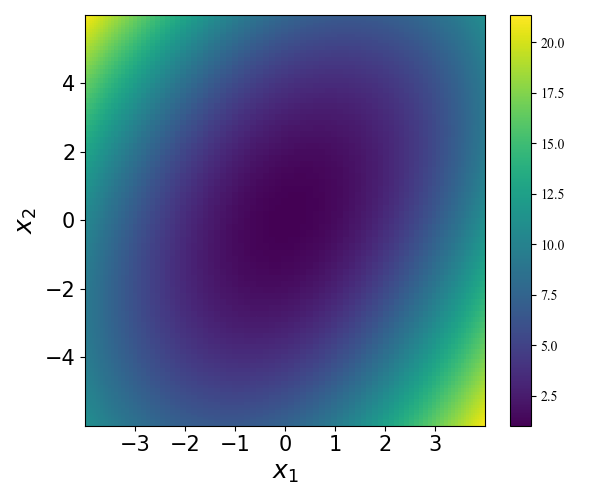}}
\hspace{0.1cm}
\subfloat[$(\bm{b}\bm{b}^\top)_{11}-(\BB_\text{NN})_{11}$]{
\includegraphics[scale=0.30]{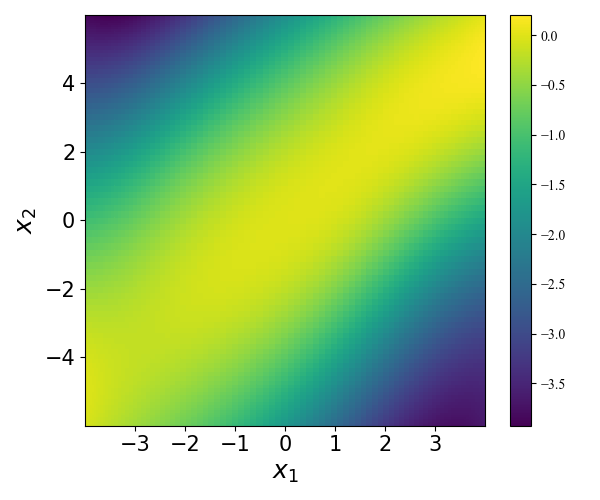}}
\caption{The comparison of first component of $\bm{b}\bm{b}^\top$. (a) $(\bm{b}\bm{b}^\top)_{11}$, (b) $(\BB_\text{NN})_{11}$ and (c) their difference.}
\label{fig:b_approximate}
\end{figure}

The relative $L_2$ error between $\Ba_\text{NN}$ and $\Ba$ is $2.94\times 10^{-2}$. Figure~\ref{fig:a_approximate} displays the spatial profile of the first components of $\Ba$,  $\Ba_\text{NN}$, and their difference on the computational domain $\Omega$. The relative $L_2$ error between $\BB_\text{NN}$ and $\Bb\Bb^\top$ is $6.98\times 10^{-2}$. To check the pointwise accuracy of the estimates, We plot the first diagonal components of $\Bb\Bb^\top$, $\BB_\text{NN}$ on the computational domain $\Omega$ and their difference in Figure~\ref{fig:b_approximate}. We can see our method works well on fitting the drift and diffusion terms.

Given the approximate drift and diffusion coefficients, we now empirically validate the result in Lemma~\ref{lemma_perterror} on the computational domain 
\BEA
\Bx^{n+1} - \Bx^n = \Ba_\text{NN}(\Bx^n)\delta t + \bm{U}(\Bx^n) \bm{S}(\Bx^n)^{\frac{1}{2}} \bm{U}(\Bx^n)^\top \sqrt{\delta t} \Bxi_n, \quad\quad \Bxi_n\sim\mathcal{N}(0,\bm{I}_2),
\label{eqn:NNeuler}
\EEA
where $\bm{U}(\Bx^n) \bm{S}(\Bx^n)\bm{U}(\Bx^n)^\top$ is the eigendecomposition of $\BB_\text{NN}(\Bx^n)$. We denote these empirical statistics to be defined with respect to the distribution $\hat{\pi}^{EM}$ that approximates $\hat{\pi}$. Compare to the ground truth statistics, the statistics
of $\hat{\pi}^{EM}$ are subjected to errors from the estimation of $\Ba$, $\BB$, and from the EM integration.
In Table~\ref{tab:EM}, we note that when $\delta t=0.05$, the covariance error of $\Ba_\text{NN},\BB_\text{NN}$ is comparable to the error of $\Ba, \Bb$. When $\delta t=0.01$, the covariance of $\Ba_\text{NN},\BB_\text{NN}$ becomes much closer to the ground truth than $\delta t=0.05$.

\begin{table}[htbp]
  \centering
  \begin{adjustbox}{width=0.8\textwidth,center}
    \begin{tabular}{ccccccc}
    \toprule
    Distribution & $\pi$ &       & $\tilde{\pi}:=\pi^\text{EM}$ &       & \multicolumn{2}{c}{$\hat{\pi}^\text{EM}$} \\
\cmidrule{1-2}\cmidrule{4-4}\cmidrule{6-7}   $\delta t$    & N/A   &       & 0.05  &       & 0.05  & 0.01 \\
    \midrule
    mean  &  $\begin{bmatrix}
0.000&0.000 \\
\end{bmatrix}$     &       &    $\begin{bmatrix}
-0.002&0.000 \\
\end{bmatrix}$   &       &    $\begin{bmatrix}
0.001&0.004 \\
\end{bmatrix}$   &  $\begin{bmatrix}
0.002&-0.003 \\
\end{bmatrix}$\\
    \midrule
    covariance &   $\begin{bmatrix}
1.000 & 0.500 \\
0.500 & 4.000
\end{bmatrix}$    &       &   $\begin{bmatrix}
1.127 & 0.499 \\
0.499 & 4.398
\end{bmatrix}$    &       &  $\begin{bmatrix}
1.115 & 0.490 \\
0.490 & 4.347
\end{bmatrix}$     &  $\begin{bmatrix}
1.013 & 0.501 \\
0.501 & 3.984
\end{bmatrix}$\\
    \bottomrule
    \end{tabular}%
    \end{adjustbox}
    \vspace{5pt}
    \caption{Comparison of mean and covariance statistics corresponding to the ground truth distribution $\pi$, discrete Markov chain induced by EM scheme in \eqref{eqn:euler}, $\tilde{\pi}:=\pi^{EM}$, and the discrete Markov chain generated by \eqref{eqn:NNeuler} for various $\delta t$ whose invariant distribution is denoted as $\hat{\pi}^\text{EM}$. }
  \label{tab:EM}%
\end{table}%

\subsubsection{Computation of the density function}
\begin{figure}
\centering
\subfloat[true density $p$]{
\includegraphics[scale=0.30]{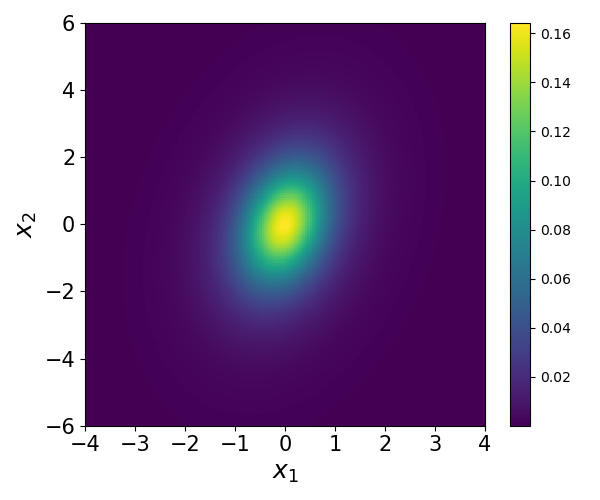}}
\hspace{0.1cm}
\subfloat[$p-\tilde{p}_\text{NN}(\cdot;\bm{\theta})$]{
\includegraphics[scale=0.30]{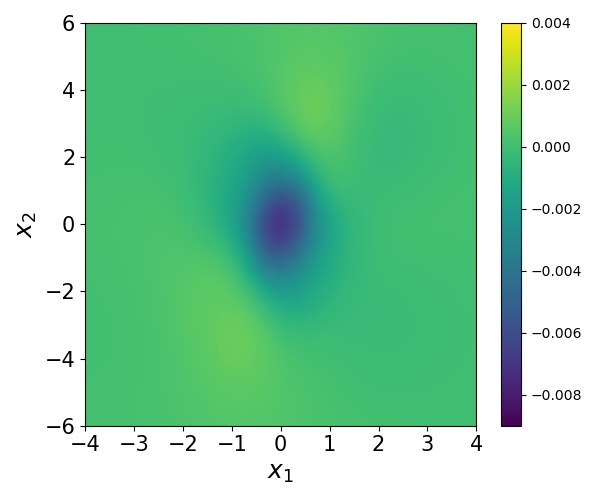}}
\hspace{0.1cm}
\subfloat[$p-\hat{p}_\text{NN}(\cdot;\bm{\theta})$]{
\includegraphics[scale=0.30]{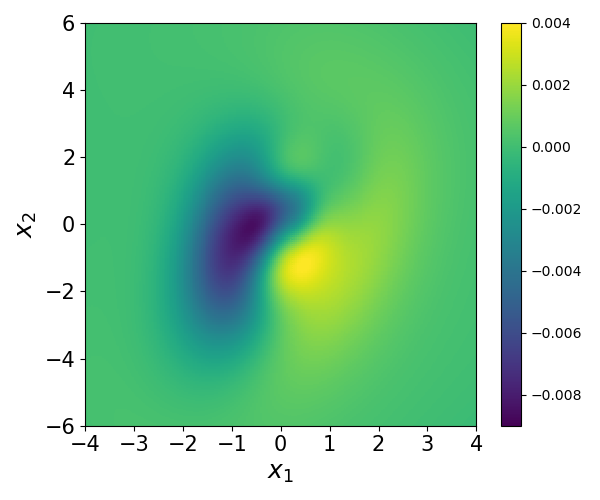}}
\caption{The comparison of solutions. (a) True density $p$, (b) difference between $p$ and $\tilde{p}_\text{NN}(\cdot;\bm{\theta})$ and (c) difference between $p$ and  $\hat{p}_\text{NN}(\cdot;\bm{\theta})$. Here $\hat{p}_\text{NN}(\cdot;\bm{\theta})$ is obtained by optimizing~\eqref{eqn:loss} with $\Ba_\text{NN}$ and $\BB_\text{NN} $, while $\tilde{p}_\text{NN}(\cdot;\bm{\theta})$ is obtained by optimizing~\eqref{eqn:loss} with $\Ba$ and $\Bb\Bb^\top$.}
\label{fig:student_solution}
\end{figure}

We optimize the loss~\eqref{eqn:loss} with $\Ba_\text{NN}$ and $\BB_\text{NN}$ and obtain the solution $\hat{p}_\text{NN}(\cdot;\bm{\theta})$. The relative $L_2$ error between $\hat{p}_\text{NN}(\cdot;\bm{\theta})$ and the true density~\eqref{student solution} is $6.78\times10^{-2}$. To quantify the error induced by the regression alone, we replace $\Ba_\text{NN}$ and $\BB_\text{NN}$ of $\hat{\mathcal{L}}^*$ in~\eqref{eqn:loss} with the underlying coefficients, $\Ba$ and $ \Bb\Bb^\top$, and optimize~\eqref{eqn:loss} with differential operator $\mathcal{L}^*$ in the first term. We denote the corresponding solution by $\tilde{p}_\text{NN}(\cdot;\bm{\theta})$. The relative $L_2$ error between $\tilde{p}_\text{NN}(\cdot;\bm{\theta})$ and the true density~\eqref{student solution} is $3.97\times10^{-2}$. We can see $\hat{p}_\text{NN}(\cdot;\bm{\theta})$ achieves the error of same magnitude as $\tilde{p}_\text{NN}(\cdot;\bm{\theta})$.  Figure~\ref{fig:student_solution} shows the true density and the differences between the true density and the network solutions $\hat{p}_\text{NN}(\cdot;\bm{\theta})$, $\tilde{p}_\text{NN}(\cdot;\bm{\theta})$, plotted as functions of the computational domain $\Omega$. Notice that the errors are more prominent when the coefficients $\Ba$ and $\Bb$ are estimated, as expected.

\subsection{The Langevin dynamics}

We consider a molecular model describing the dynamics of $M$ atoms with mass 1. We assume the $M$ particles are spaced in a chain with a periodic boundary condition. Let the equilibrium distance between two neighboring particles be $a_0$, then the equilibrium position of the $m$-th particle is $ma_0$. Denote $r_m$ as the displacement of the $m$-th particle from its equilibrium position, and denote $v_m$ as its velocity. The Langevin dynamics of this model is described as follows
\begin{equation}\label{31}
\begin{split}
\dot{\Bv}&=-\nabla_{\Br} U(\Br)-\gamma\Bv+\sqrt{2k_BT\gamma}\dot{\BW}_t,\\
\dot{\Br}&=\Bv,
\end{split}
\end{equation}
where $\Bv=[v_1,\cdots,v_M]^\top$ and $\Br=[r_1,\cdots,r_M]^\top$ are the velocities and displacement of all particles; $\BW_t=[W_t^{(1)},\cdots,W_t^{(M)}]^\top$ is an $M$-dimensional Wiener process; $U$ is some potential function; $\gamma$ is the friction constant; $k_BT$ is the temperature. The mass of particles is set to be unity in \eqref{31}. The equilibrium distribution of \eqref{31} is given by
\begin{equation}\label{36}
p(\Bv,\Br)\propto\exp\left[-\frac{1}{k_BT}\left(U(\Br)+\frac{1}{2}|\Bv|^2\right)\right].
\end{equation}

In the numerical simulation, we take the Lennard-Jones potential \cite{Ishimori1982}, which is given by
\begin{equation}\label{37}
U(\Br)=\sum_{i=1}^{M}\sum_{j=i-2}^{i-1}\psi(r_i-r_j+(i-j)a_0),\quad r_0:=r_M,~r_{-1}:=r_{M-1}
\end{equation}
with
\begin{equation}
\psi(r)=|r|^{-12}-|r|^{-6}.
\end{equation}
The model parameters of this example is set to be $a_0=1$, $\gamma=0.5$, $k_BT=0.25$, $M=10$.

\subsubsection{Data generation}
We generate the data by Euler-Maruyama discretization, namely,
\begin{equation}
\begin{split}\label{EMLangevin}
\Bv^{n+1}&=\Bv^n-(\nabla_{\Br^n} U(\Br^n)+\gamma\Bv^n)\delta t+\sqrt{2k_BT\gamma\delta t}\Bxi,\quad\Bxi\sim\left(\mathcal{N}(0,1)\right)^M,\\
\Br^{n+1}&=\Br^n+\Bv^n\delta t,
\end{split}
\end{equation}
for $n=0,1,\cdots,N-1$ with the initial states
\begin{equation*}
\Bv^0=0,\quad\Br^0\sim\left(\mathcal{N}(0,0.01)\right)^M.
\end{equation*}

\begin{figure}
\centering
\subfloat[distribution of $(r_1,r_2)$]{
\includegraphics[scale=0.5]{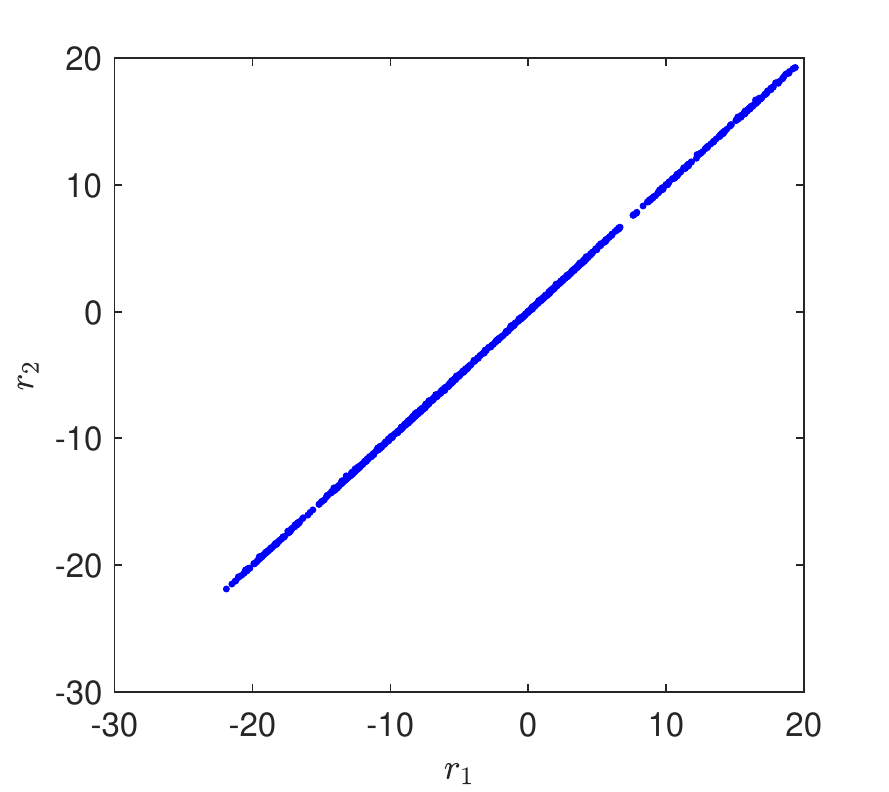}}
\subfloat[distribution of $(d_1,d_2)$]{
\includegraphics[scale=0.5]{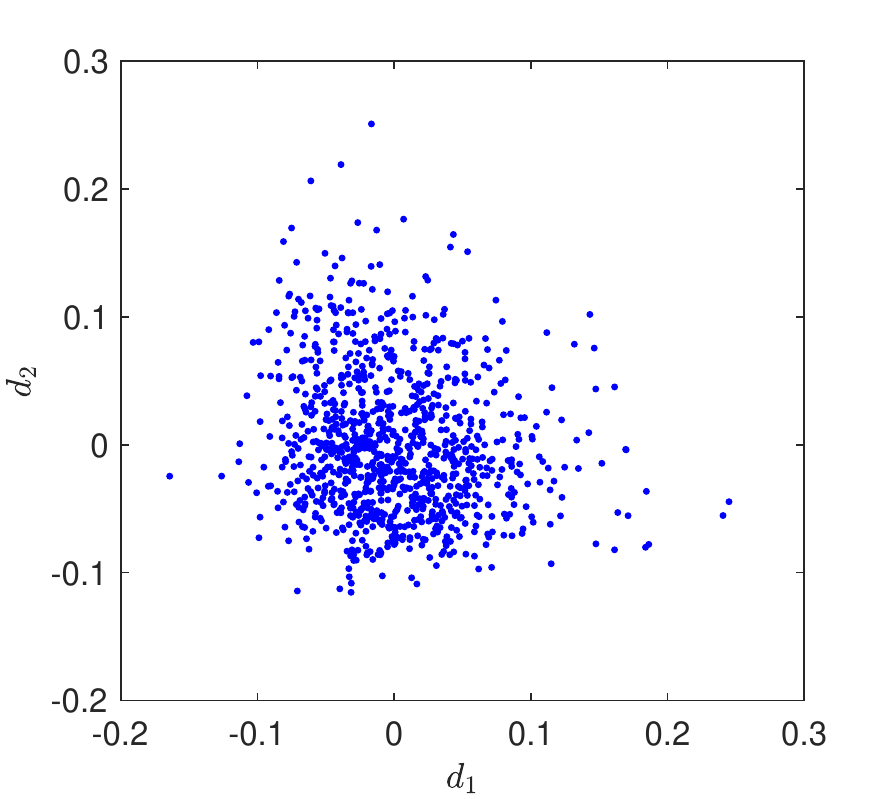}}
\caption{\em The distribution of the original dataset $\mathcal{X}$ in the $(r_1,r_2)$-plane and the distribution of the transformed dataset $\hat{\mathcal{X}}$ in the $(d_1,d_2)$-plane.}
\label{Fig_Lang_data_plot}
\end{figure}

In this example, we set $\delta t=0.0005$ and $N=10^7$. Following the notation in Section~\ref{sec:regression}, we denote $\mathcal{X}:=\{\Bv^n,\Br^n\}_{n=0}^N$ as the original data set. If we visualize the distribution of $\mathcal{X}$ by projecting it onto the $(r_1,r_2)$-plane (Figure \ref{Fig_Lang_data_plot}), it is observed that displacement components are distributed near a straight line.
To simplify the computation and visualization, we consider a coordinate transformation that will map $\mathcal{X}$ to be enclosed by a hyper-rectangle. Specifically, we introduce the following coordinate transformation,
\begin{equation*}
\mathcal{T}:\mathbb{R}^M\rightarrow\mathbb{R}^{M-1},\quad\Bd:=[d_1,\cdots,d_{M-1}]^\top=\mathcal{T}(\Br)=\left[r_2-r_1,\cdots,r_M-r_{M-1}\right]^\top,
\end{equation*}
where $\Bd$ is called the relative displacement. Note that the map $\mathcal{T}$ implies $r_1-r_M=-\sum_{m=1}^{M-1}d_m$. If we define the transformed dateset  $\hat{\mathcal{X}}:=\{\Bv^n,\Bd^n\}_{n=0}^N$ with $\Bd^n=\mathcal{T}(\Br^n)$ and project it onto the $(d_1,d_2)$-plane (Figure \ref{Fig_Lang_data_plot}), then it is observed that most points in $\hat{\mathcal{X}}$ are located near the origin and form a circular region. Consequently, we apply the proposed method to the transformed dataset $\hat{\mathcal{X}}$ in the practical computation.

\subsubsection{Identification of the drift and diffusion coefficients.}
Now we aim to identify the drift term $\Ba(\Bv,\Br)$ and the diffusion $\Bb\Bb^\top$ of the underlying dynamics. Due to transformation $\mathcal{T}$, we define $\hat{\Ba}(\Bv,\Bd):=\Ba(\Bv,\Br)$ and aim to identify $\hat{\Ba}$ by the optimization \eqref{eqn:traininga} using the dataset $\hat{\mathcal{X}}$. Note $\hat{\Ba}(\Bv,\Bd)$ is a vector-valued function with $(2M-1)$-dimensional inputs and $2M$-dimensional outputs. In this example, to obtain higher accuracy, we use an individual neural network with $(2M-1)$-dimensional inputs and scalar outputs to approximate the each component of $\hat{\Ba}(\Bv,\Bd)$, solving the regression problem in \eqref{eqn:traininga}. In this application, this is a regression over training data set $(\hat{\mathcal{X}},\mathcal{Y})$, where $\mathcal{Y}:=\{\frac{\Bv^{n+1}-\Bv^n}{\delta t},\frac{\Br^{n+1}Br^n}{\delta t}\}_{n=0}^{N-1}$

\comment{
, denoted as $\hat{a}_\text{NN}(\Bv,\Bd;\Btheta_k)$, to approximate the $k$-th component of $\hat{\Ba}(\Bv,\Bd)$, say $\hat{a}_k(\Bv,\Bd)$, where $k=1,\cdots,2M$. Specifically, we train the networks $\{\hat{a}_\text{NN}(\Bv,\Bd;\Btheta_k)\}$ by the following optimizations,
\begin{gather}
\underset{\Btheta_k}{\min}\overset{N-1}{\underset{n=0}{\sum}}\left|\frac{v_k^{n+1}-v_k^n}{\delta t}-\hat{a}_\text{NN}(\Bv^n,\Bd^n;\Btheta_k)\right|^2,\quad\text{for~}k=1,\cdots,M,\label{34_1}\\
\underset{\Btheta_k}{\min}\overset{N-1}{\underset{n=0}{\sum}}\left|\frac{r_{k-M}^{n+1}-r_{k-M}^n}{\delta t}-\hat{a}_\text{NN}(\Bv^n,\Bd^n;\Btheta_k)\right|^2,\quad\text{for~}k=M+1,\cdots,2M.\label{34_2}
\end{gather}}

In practice. we set each component of $\hat{\Ba}_\text{NN}$ to be a fully connected ReLU network with 3 layers and 100 neurons in each layer. We employ Adams optimizer with 1000 epochs, and the learning rates are set to decay from $10^{-3}$ to $10^{-5}$. The relative $\ell^2$ training errors for the first $M$ components corresponding to the velocity are observed to be between $3.87\times10^{-2}$ and $5.60\times10^{-2}$, and the errors for the next $M$ components are between $5.04\times10^{-5}$ and $8.99\times10^{-5}$.

\comment{
\begin{table}[htbp]
  \centering
    \begin{tabular}{cccccc}
    \toprule
    $k$ & 1 & 2 & 3 & 4 & 5\\\midrule
    relative $\ell_2$ error & $4.28\times10^{-2}$ & $5.60\times10^{-2}$ & $3.98\times10^{-2}$ & $4.27\times10^{-2}$ & $3.70\times10^{-2}$ \\\midrule
    $k$ & 6 & 7 & 8 & 9 & 10\\\midrule
    relative $\ell_2$ error & $4.36\times10^{-2}$ & $4.76\times10^{-2}$ & $4.32\times10^{-2}$ & $4.14\times10^{-2}$ & $3.87\times10^{-2}$ \\\midrule
    $k$ & 11 & 12 & 13 & 14 & 15\\\midrule
    relative $\ell_2$ error & $6.46\times10^{-5}$ & $5.85\times10^{-5}$ & $5.04\times10^{-5}$ & $7.37\times10^{-5}$ & $5.32\times10^{-5}$ \\\midrule
    $k$ & 16 & 17 & 18 & 19 & 20\\\midrule
    relative $\ell_2$ error & $5.61\times10^{-5}$ & $8.99\times10^{-5}$ & $5.75\times10^{-5}$ & $6.58\times10^{-5}$ & $6.04\times10^{-5}$ \\
    \bottomrule
    \end{tabular}%
    \caption{$\ell^2$ errors between the identified drift function $\hat{a}_\text{NN}(\Bv',\Bd';\Btheta_k)$ and the true drift function $\hat{\Ba}_k(\Bv,\Bd)$ in the Langevin dynamics example.}
    \label{tab:a_error_Langevin}
\end{table}
}

Next, we consider the approximation $\BB_{NN}$ to the constant matrix $\Bb\Bb^\top$ using the formula in \eqref{constantBBtop}. In this example, since $\Bb\Bb^\top$ is a diagonal matrix, we also set $\BB_\text{NN}$ to be diagonal with components $(b_{11},\ldots,b_{2M,2M})$. The errors $|b_{kk}-(\Bb\Bb^\top)_{kk}|$ for the first $M$ components are observed to be between $6.32\times10^{-6}$ and $2.67\times10^{-6}$, and the errors for the next $M$ components are between $8.07\times10^{-13}$ and $3.63\times10^{-12}$.

\comment{
\begin{table}[htbp]
  \centering
    \begin{tabular}{cccccc}
    \toprule
    $k$ & 1 & 2 & 3 & 4 & 5\\\midrule
    error & $7.95\times10^{-5}$ & $1.67\times10^{-4}$ & $2.30\times10^{-4}$ & $2.67\times10^{-4}$ & $4.98\times10^{-5}$ \\\midrule
    $k$ & 6 & 7 & 8 & 9 & 10\\\midrule
    error & $9.29\times10^{-5}$ & $2.98\times10^{-5}$ & $1.45\times10^{-4}$ & $6.32\times10^{-6}$ & $2.46\times10^{-4}$ \\\midrule
    $k$ & 11 & 12 & 13 & 14 & 15\\\midrule
    error & $1.58\times10^{-12}$ & $1.84\times10^{-12}$ & $1.16\times10^{-12}$ & $1.35\times10^{-12}$ & $8.07\times10^{-13}$ \\\midrule
    $k$ & 16 & 17 & 18 & 19 & 20\\\midrule
    error & $1.23\times10^{-12}$ & $3.63\times10^{-12}$ & $1.34\times10^{-12}$ & $1.97\times10^{-12}$ & $1.30\times10^{-12}$ \\
    \bottomrule
    \end{tabular}%
    \caption{Errors $|b_{kk}-(\Bb\Bb^\top)_{kk}|$ between the identified diffusion value $b_{kk}$ and the true diffusion value $(\Bb\Bb^\top)_{kk}$ in the Langevin dynamics example.}
    \label{tab:b_error_Langevin}
\end{table}
}

Similar to the previous example, we simulate the dynamics by the obtained $\hat{\Ba}_\tNN$ and $\BB_\tNN$,
\BEA\label{eqn:NNeuler2}
\begin{aligned}
\Bv^{n+1} - \Bv^n &= (\hat{\Ba}_\text{NN})_{1:M}(\Bv^n,\Bd^n)\delta t + (\BB_\tNN)^{\frac{1}{2}} \sqrt{\delta t} \Bxi_n, \quad\quad \Bxi_n\sim\mathcal{N}(0,\bm{I}_M),\\
\Br^{n+1} - \Br^n &= (\hat{\Ba}_\text{NN})_{M+1:2M}(\Bv^n,\Bd^n)\delta t,
\end{aligned}
\EEA
and compare it with the ground truth. For the covariance of $\pi$, Monte Carlo integration with $10^8$ points is used. For the statistics of $\tilde{\pi}$ and $\hat{\pi}^\text{EM}$, we generate a sequence of $10^7$ points. The information is shown in Table \ref{tab:EM2} for the components $\Bv_1$ and $\Bd_1$. Notice that in this case, the statistical error for estimating $\hat{\pi}^{EM}$ is not much worse than the Monte-Carlo error of $\tilde{\pi}$.


\begin{table}[htbp]
  \centering
  \begin{adjustbox}{width=0.85\textwidth,center}

   \begin{tabular}{ccccccc}
    \toprule
    Distribution &       & $\pi$ &       & $\tilde{\pi}$ &       & $\hat{\pi}^\text{EM}$ \\
\cmidrule{1-1}\cmidrule{3-3}\cmidrule{5-5}\cmidrule{7-7}    $\delta t$ &       & N/A &       & 0.0005 &       & 0.0005 \\
\cmidrule{1-1}\cmidrule{3-3}\cmidrule{5-5}\cmidrule{7-7}    mean &       & $\begin{bmatrix}0&0 \\\end{bmatrix}$ &       & $\begin{bmatrix}-0.00363&-0.00013 \\\end{bmatrix}$ &       & $\begin{bmatrix}-0.00153&-0.00003 \\\end{bmatrix}$ \\
\cmidrule{1-1}\cmidrule{3-3}\cmidrule{5-5}\cmidrule{7-7}    covariance &       & $\begin{bmatrix}0.40229 & -0.01749 \\-0.01749 & 0.00245 \end{bmatrix}$ &       & $\begin{bmatrix}0.37816 & 0.00008 \\0.00008 & 0.00292\end{bmatrix}$ &       & $\begin{bmatrix}0.40916 & 0.00041 \\0.00041 & 0.00314\end{bmatrix}$ \\
    \bottomrule
    \end{tabular}%

    \end{adjustbox}
    \caption{Comparison of mean and covariance statistics ($\Bv_1$ and $\Bd_1$) corresponding to the ground truth distribution $\pi$ , discrete Markov chain induced by EM scheme in \eqref{eqn:euler}, $\tilde{\pi}$, and the discrete Markov chain generated by \eqref{eqn:NNeuler2} for $\delta t=0.0005$ whose invariant distribution is denoted as $\hat{\pi}^\text{EM}$. }
    \label{tab:EM2}.
\end{table}%

\subsubsection{Computation of the density function}
In this section, we aim to recover the equilibrium density function based on the obtained $\{\hat{a}_\text{NN}(\Bv,\Bd;\Btheta_k)\}$ and $\BB_\text{NN}$. We let $p(\Bv,\Br)$ be the original density function in $(\Bv,\Br)$-coordinates, and define $\hat{p}(\Bv,\Bd):=p(\Bv,\Br)$ be the density function under transformation $\mathcal{T}$. Since $p(\Bv,\Br)$ satisfies \eqref{02}, we can derive the PDE for $\hat{p}(\Bv,\Bd)$, which is given by
\begin{multline}\label{32}
-\underset{k=1}{\overset{M}{\sum}}\frac{\partial}{\partial v_k}(\hat{p}\hat{a}_k)-\underset{k=M+1}{\overset{2M-1}{\sum}}\frac{\partial}{\partial d_{k-M}}\left(\hat{p}(\hat{a}_{k+1}-\hat{a}_k)\right)\\
+\frac{1}{2}\underset{k=1}{\overset{M}{\sum}}(\Bb\Bb^\top)_{kk}\frac{\partial^2}{\partial v_k^2}\hat{p}+\frac{1}{2}(\Bb\Bb^\top)_{M+1,M+1}\frac{\partial^2}{\partial d_1^2}\hat{p}\\
+\frac{1}{2}\underset{k=2}{\overset{M-1}{\sum}}(\Bb\Bb^\top)_{k+M,k+M}\left(\frac{\partial}{\partial d_k}-\frac{\partial}{\partial d_{k-1}}\right)^2\hat{p}+\frac{1}{2}(\Bb\Bb^\top)_{2M,2M}\frac{\partial^2}{\partial d_{M-1}^2}\hat{p}=0,
\end{multline}
where $\hat{a}_k$ denotes the $k$-th component of $\hat{\Ba}$.

Once the drift and diffusion coefficients are estimated, we substitute $\hat{a}_k$ with the $k$th FNN estimate, denoted as $\hat{a}_\text{NN}(\Bv,\Bd;\Btheta_k)$, and $(\Bb\Bb^\top)_{k,k}$ with the diagonal components of the estimated diffusion matrix, $b_{kk}:=(\BB_\text{NN})_{kk}$, such that \eqref{32}, becomes,
 \begin{multline}\label{33}
-\underset{k=1}{\overset{M}{\sum}}\frac{\partial}{\partial v_k}(\hat{p}\hat{a}_\text{NN}(\Bv,\Bd;\Btheta_k))-\underset{k=M+1}{\overset{2M-1}{\sum}}\frac{\partial}{\partial d_{k-M}}\left(\hat{p}(\hat{a}_\text{NN}(\Bv,\Bd;\Btheta_{k+1})-\hat{a}_\text{NN}(\Bv,\Bd;\Btheta_k))\right)\\
+\frac{1}{2}\left(\underset{k=1}{\overset{M}{\sum}}b_{kk}\frac{\partial^2\hat{p}}{\partial v_k^2}+b_{M+1,M+1}\frac{\partial^2\hat{p}}{\partial d_1^2}
+\underset{k=2}{\overset{M-1}{\sum}}b_{k+M,k+M}\left(\frac{\partial}{\partial d_k}-\frac{\partial}{\partial d_{k-1}}\right)^2\hat{p}+b_{2M,2M}\frac{\partial^2\hat{p}}{\partial d_{M-1}^2}\right)=0,
\end{multline}

\begin{figure}
\centering
\includegraphics[scale=0.5]{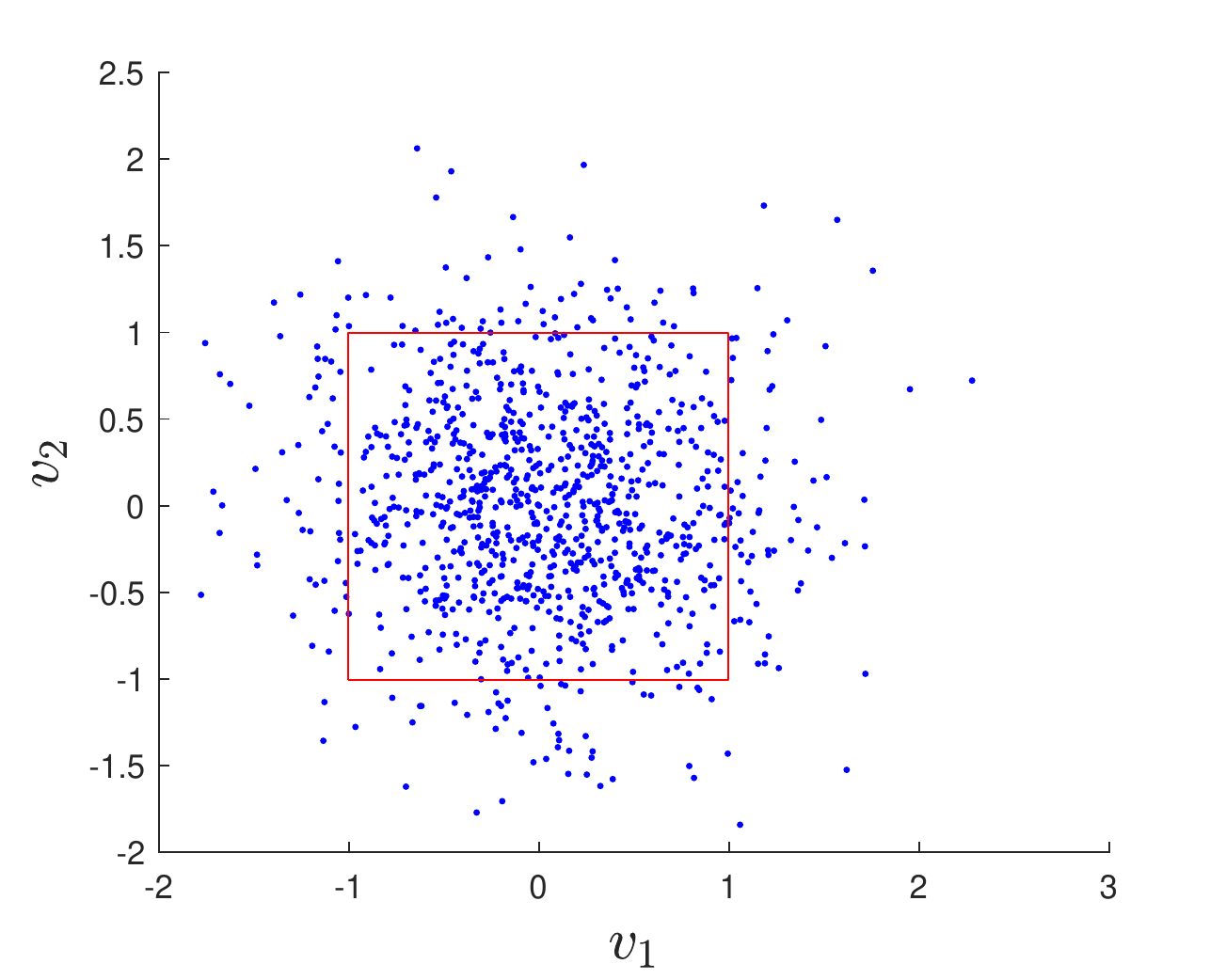}
\includegraphics[scale=0.5]{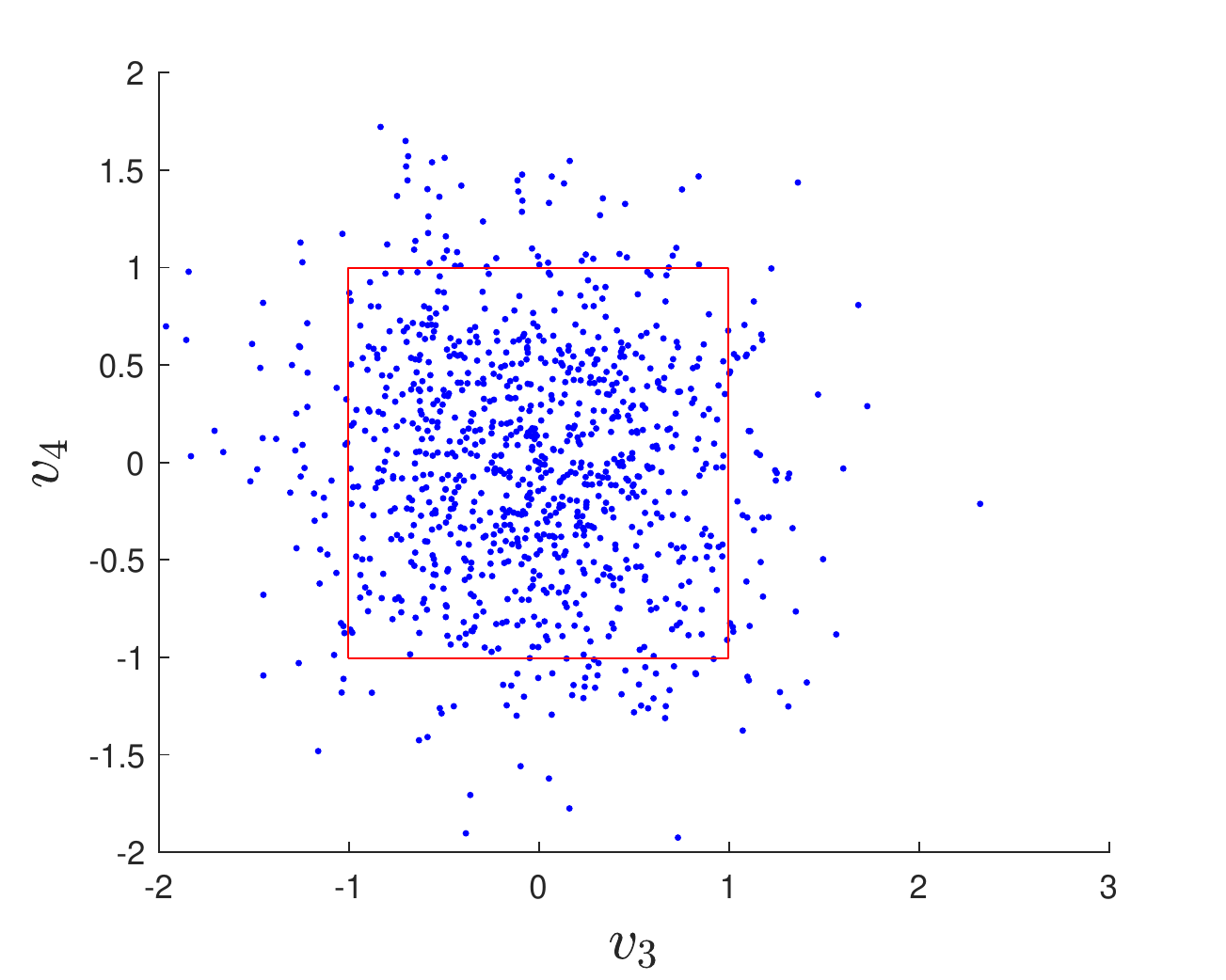}\\
\includegraphics[scale=0.5]{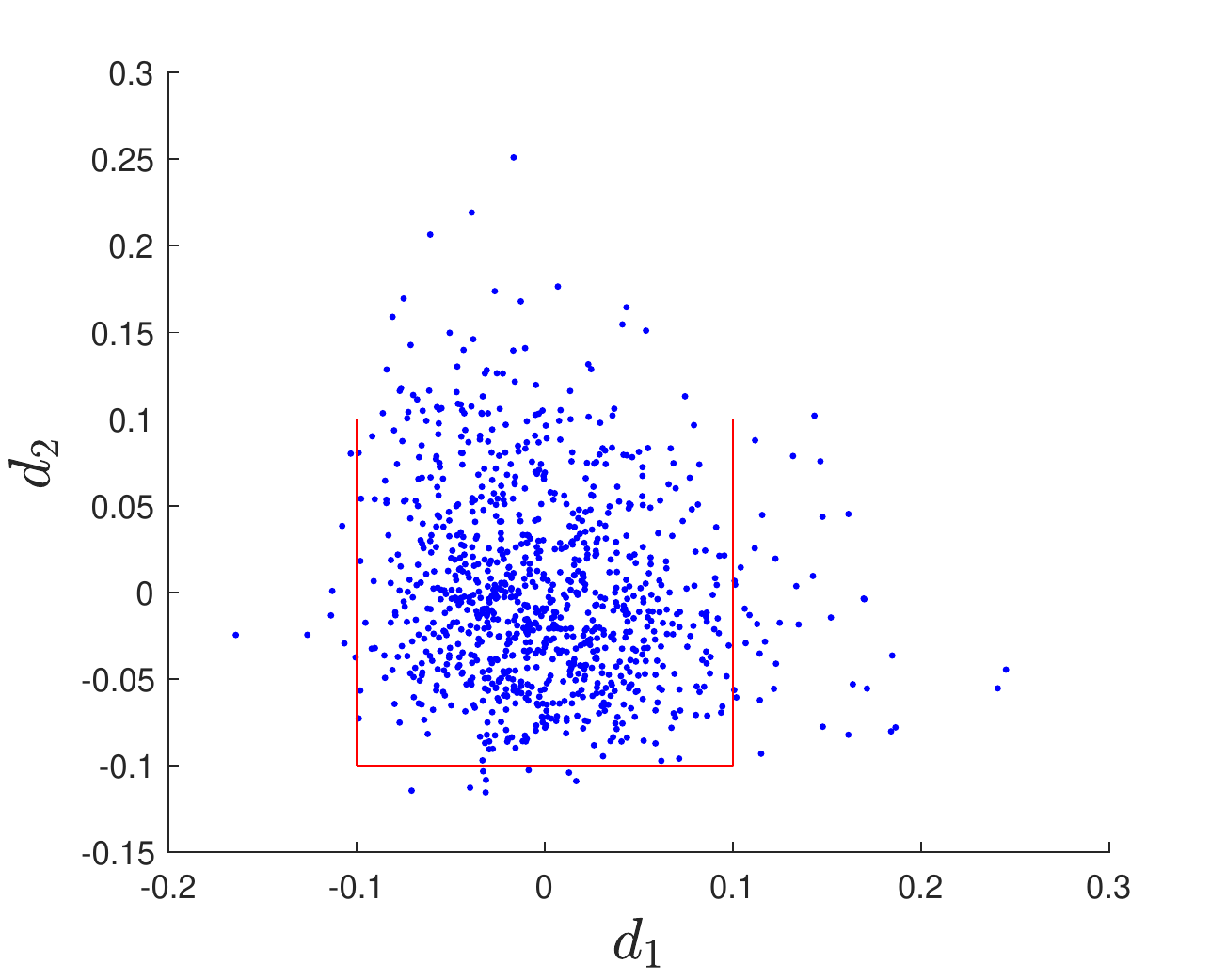}
\includegraphics[scale=0.5]{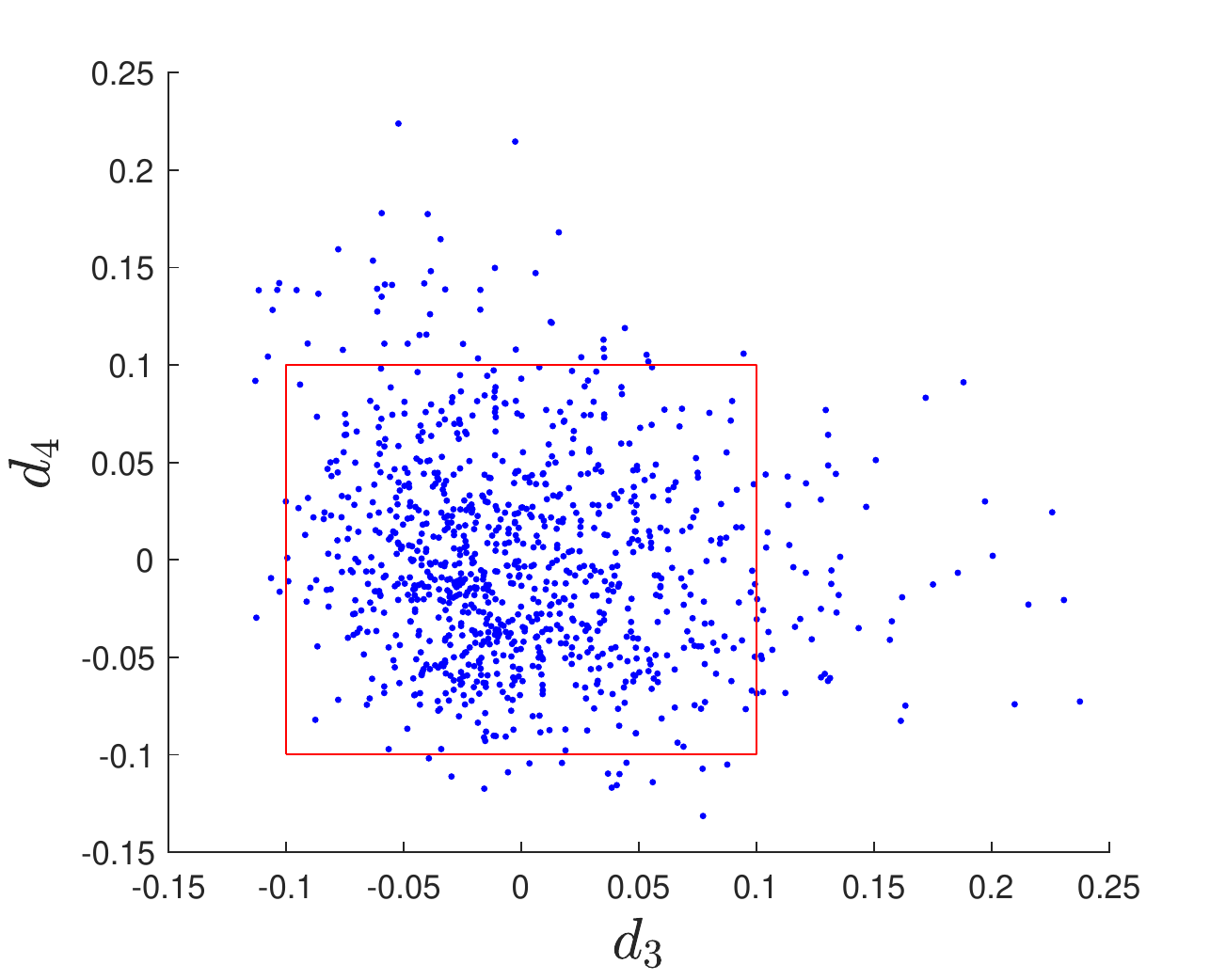}
\caption{\em The projections of the dataset $\mathcal{X}$ (blue points) and the enclosing region $\Omega$ (red boxes) onto $(v_1,v_2)$, $(v_3,v_3)$, $(d_1,d_2)$, $(d_3,d_4)$-planes}
\label{Fig_Lang_data_region_plot}
\end{figure}

Next, we select a bounded domain in which the PDE \eqref{33} will be solved. Our choice is to use a hyperrectangle $\Omega=\underset{k=1}{\overset{2M-1}{\prod}}[c_k-s_k,c_k+s_k]$ to enclose most of the points in $\hat{\mathcal{X}}$. At the same time, we expect $\Omega$ to be also densely covered by the points in $\hat{\mathcal{X}}$. By this principle, we set $c_k$ as the component-wise mean of the points in $\hat{\mathcal{X}}$, namely,
\begin{equation}
c_k=\begin{cases}
\frac{1}{N}\sum_{n=1}^{N} v_k,\quad\text{for~}k=1,\cdots,M,\\
\frac{1}{N}\sum_{n=1}^{N} d_k,\quad\text{for~}k=M+1,\cdots,2M-1,
\end{cases}
\end{equation}
and set $s_k$ empirically as follows.
\begin{equation}
s_k=\begin{cases}
1.0,\quad\text{for~}k=1,\cdots,M,\\
0.1,\quad\text{for~}k=M+1,\cdots,2M-1.
\end{cases}
\end{equation}
For clarity, we display the projections of $\hat{\mathcal{X}}$ and $\Omega$ onto coordinate planes in Figure \ref{Fig_Lang_data_region_plot}.

\begin{figure}
\centering
\includegraphics[scale=0.8]{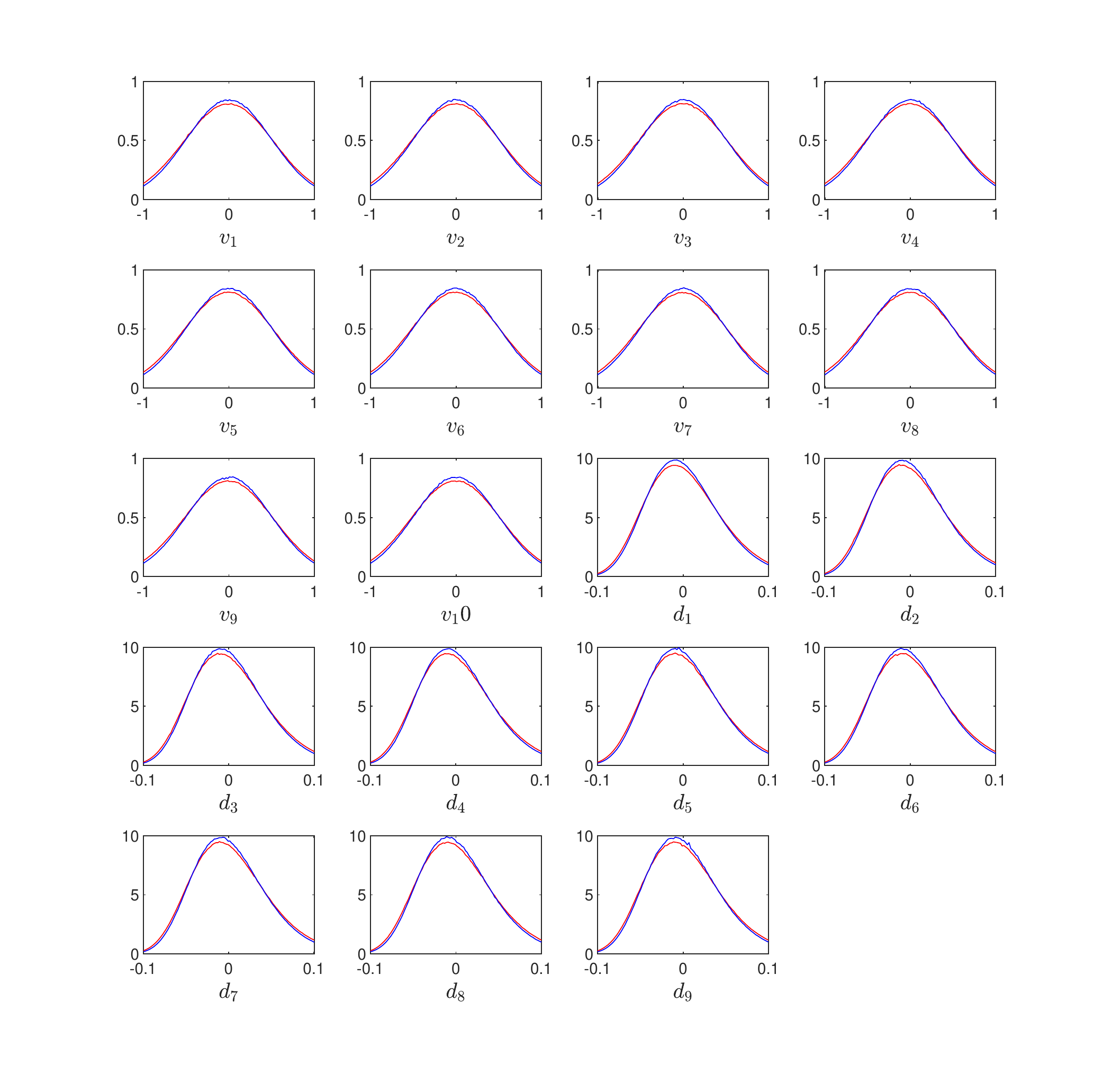}
\caption{\em Marginal densities of the computed density function $\hat{p}_\text{NN}$ (red curves) and the true density function $\hat{p}$ (blue curves) for all components.}
\label{Fig_Lang_marginal_density}
\end{figure}

We take a neural network $\hat{p}_\text{NN}(\Bv,\Bd;\Btheta)$ to approximate $\hat{p}(\Bv,\Bd)$. Then we solve the PDE \eqref{33} with the least square method introduced in Section \ref{Sec_least_square} to determine $\hat{p}_\text{NN}(\Bv,\Bd;\Btheta)$.
Specifically, we solve the least-square problem in \eqref{pdeempiricalloss} with $\gamma=0$, ignoring the artificial boundary constraint since the function values at the boundary $\partial\Omega$ are small, they range from $7\times 10^{-7}$ to $4\times 10^{-6}$. Meanwhile, $90\%$ of the points in $\hat{\mathcal{X}}\cap\Omega$ are selected as the training set, denoted as $\hat{D}_\text{T}$, and the other $10\%$ are chosen as the testing set, denoted as $\hat{D}_\text{S}$, for the evaluation of the solution error. In practice. we set each $\hat{p}_\text{NN}$ to be a fully-connected network having 3 layers and 100 neurons in each layer with activation function $\max\{x^3,0\}$. Adams optimizer is used to solve the optimization with 1000 epochs, and the learning rates are set to decay from $10^{-4}$ to $10^{-5}$. Once $\hat{p}_\text{NN}$ is obtained, we evaluate the result by computing the error between $\hat{p}_\text{NN}(\Bv,\Bd)$ and the true density function $\hat{p}(\Bv,\Bd)$. From \eqref{36}-\eqref{37}, we directly have the expression of $\hat{p}(\Bv,\Bd)$, namely,
\begin{equation}
\hat{p}(\Bv,\Bd)=c\cdot\hat{p}_0(\Bv,\Bd):=c\cdot\exp\left[-\frac{1}{k_BT}\left(\hat{U}(\Bd)+\frac{1}{2}|\Bv|^2\right)\right]
\end{equation}
with
\begin{multline}
\hat{U}(\Bd)=\psi(-\sum_{i=1}^{M-1}d_i+a_0)+\psi(-\sum_{i=1}^{M-2}d_i+2a_0)+\psi(d_1+a_0)+\psi(-\sum_{i=2}^{M-1}d_i+2a_0)\\
+\sum_{i=3}^{M}\psi(d_{i-1}+a_0)+\sum_{i=3}^{M}\psi(d_{i-1}+d_{i-2}+2a_0),
\end{multline}
where $c$ is determined by the uniform condition
\begin{equation}\label{38}
c=\left(\int_{\mathbb{R}^{2M-1}}\hat{p}_0(\Bv,\Bd)\right)^{-1}.
\end{equation}
Then the relative $\ell^2$ error between $\hat{p}_\text{NN}(\Bv,\Bd)$ and $\hat{p}(\Bv,\Bd)$ can be computed according to \eqref{relativeL2error} with $L^2(\Omega)$ replaced by $L^2(\hat{D}_S)$, where the integral is replaced by an average over the testing data set $\hat{D}_S$. In this numerical result, we found thatthe relative $\ell^2$ error of the computed density function $\hat{p}_\text{NN}$ is $5.402\times10^{-2}$. In Figure \ref{Fig_Lang_marginal_density}, we also show
\comment{can be computed by
\begin{equation}
\frac{\underset{(\Bv',\Bd')\in\hat{D}_\text{S}}{\sum}\left|\hat{p}(\Bv',\Bd')-\hat{p}_\text{NN}(\Bv',\Bd';\Btheta)\right|^2}{\underset{(\Bv',\Bd')\in\hat{D}_\text{S}}{\sum}\left|\hat{p}(\Bv',\Bd')\right|^2}.
\end{equation}
Finally, the relative $\ell^2$ error of the computed density function $\hat{p}_\text{NN}$ is $5.402\times10^{-2}$.  We also show in Figure \ref{Fig_Lang_marginal_density}
}
the marginal densities of $\hat{p}_\text{NN}$
\begin{equation}\label{39}
\begin{split}
\hat{p}^{\text{marginal}}_{\text{NN},k}(v_k)&:=\int_{(\Bv,\Bd)\backslash v_k\in\mathbb{R}^{2M-2}}\hat{p}_\text{NN}(\Bv,\Bd;\Btheta),\quad\text{for~}k=1,\cdots,M,\\
\hat{p}^{\text{marginal}}_{\text{NN},k}(d_k)&:=\int_{(\Bv,\Bd)\backslash d_k\in\mathbb{R}^{2M-2}}\hat{p}_\text{NN}(\Bv,\Bd;\Btheta),\quad\text{for~}k=M+1,\cdots,2M-1,
\end{split}
\end{equation}
compared with the following true marginal densities
\begin{equation}\label{40}
\begin{split}
\hat{p}^{\text{marginal}}_k(v_k)&:=\int_{(\Bv,\Bd)\backslash v_k\in\mathbb{R}^{2M-2}}\hat{p}(\Bv,\Bd),\quad\text{for~}k=1,\cdots,M,\\
\hat{p}^{\text{marginal}}_k(d_k)&:=\int_{(\Bv,\Bd)\backslash d_k\in\mathbb{R}^{2M-2}}\hat{p}(\Bv,\Bd),\quad\text{for~}k=M+1,\cdots,2M-1,
\end{split}
\end{equation}
where the integrals in \eqref{38}, \eqref{39} and \eqref{40} are computed by the Monte Carlo method. Notice the accurate estimation of the marginal densities of the velocity components that are Gaussian and the marginal densities of the relative displacement components that are non-symmetric.

\section{Conclusion}\label{sec:cond}
In this paper, we developed a deep learning-based method to estimate the stationary density of an unknown It\^o diffusion SDE from a time series induced by the Euler-Maruyama solver. Neural networks are employed to approximate the drift, diffusion, and stationary density of the underlying dynamics. In our method, the first step is learning the drift and diffusion coefficients by solving least square regressions corresponding to the available dataset, and the second step is solving the steady-state Fokker-Plank equation formed by the estimated drift and diffusion coefficients. Theoretically, we deduced an error bound for the proposed approach for an SDE with global Lipschitz drift coefficients and constant diffusion matrix, accounting errors contributed by the discretization of the SDE in the training data, the regression of the drift terms using fully-connected ReLU networks with arbitrary width and layers, and the regression solution to the Fokker-Planck PDE using a fully-connected two-layer neural network with the $\text{ReLU}^3$ activation function. This error bound is deduced under various assumptions that underpin the perturbation theory result in \cite{zhang2021error}, generalization errors in approximating Lipschitz continuous functions in \cite{Jiao2021} and in solving PDEs in \cite{LuoYang2020}.

From this theoretical study, we observe two difficult aspects that warrant careful treatments in future studies. The first issue is concerning the incompatibility of the topologies that characterize the perturbation theory and machine learning generalization theory. Since the bound in \eqref{perturbbound} is stronger than an $L^2$-error bound in generalization theory, one requires a tacit assumption of consistency in the sense of \eqref{consistentestimator}, which is not easily verified in practice. The second issue is concerning the incompatibility of the computational and physical domains, which is admitted under the Assumption\ref{assumpdomain}. Particularly, while the underlying stochastic process is defined on $\mathbb{R}^d$, the error estimation that accounts for finite samples the training for $\Ba$ and $\hat{p}$ is not easily guaranteed for the entire unbounded domain. Besides, it is also only feasible to employ the computation over a bounded domain.

In numerical simulations, we verified the effectiveness of the proposed method on two examples, a two-dimensional Student’s t-distribution, and the 20-dimensional Langevin dynamics. Although the proposed data-driven methods show encouraging numerical results on the approximation of the invariant statistics and densities, the empirical loss function in \eqref{pdeempiricalloss} requires samples $\Bx_{I}, \Bx_{II}^n$ and $\Bx_{III}^n$. Such a requirement may not be viable when the geometry is more complicated than hypercubes. While sampling the first term in \eqref{pdeempiricalloss} is avoidable by a Monte-Carlo over the available time series as we have done in our numerical examples, generating samples for the second and third terms in the loss function in \eqref{pdeempiricalloss} is unavoidable. In the future, we plan to consider different penalties such as the one proposed in \cite{zhai2020deep} which requires no additional samples other than the available time series.



\section*{Acknowledgment}

The research of JH was partially supported under the NSF grant DMS-1854299. HY was partially supported by the US National Science Foundation under award DMS-1945029.

\bibliographystyle{plain}

\begin{thebibliography}{10}

\bibitem{averina1988numerical}
TA~Averina and SS~Artemiev.
\newblock Numerical solution of systems of stochastic differential equations.
\newblock {\em Russian Journal of Numerical Analysis and Mathematical
  Modelling}, 3(4):267--286, 1988.

\bibitem{barron1993}
A.~R. Barron.
\newblock Universal approximation bounds for superpositions of a sigmoidal
  function.
\newblock {\em IEEE Transactions on Information Theory}, 39(3):930--945, May
  1993.

\bibitem{Caragea2020}
A.~Caragea, P.~Petersen, and F.~Voigtlaender.
\newblock {Neural network approximation and estimation of classifiers with
  classification boundary in a Barron class}.
\newblock {\em arXiv e-prints}, arXiv:2011.09363, 2020.

\bibitem{chavanis2008nonlinear}
Pierre-Henri Chavanis.
\newblock Nonlinear mean field fokker-planck equations. application to the
  chemotaxis of biological populations.
\newblock {\em The European Physical Journal B}, 62(2):179--208, 2008.

\bibitem{Chen2019NonparametricRO}
Minshuo Chen, Haoming Jiang, Wenjing Liao, and T.~Zhao.
\newblock Nonparametric regression on low-dimensional manifolds using deep relu
  networks.
\newblock {\em arXiv: Learning}, 2019.

\bibitem{chen2021representation}
Ziang Chen, Jianfeng Lu, and Yulong Lu.
\newblock On the representation of solutions to elliptic pdes in barron spaces.
\newblock {\em arxiv:2106.07539}, 2021.

\bibitem{Chen1}
Zixiang Chen, Yuan Cao, Difan Zou, and Quanquan Gu.
\newblock How much over-parameterization is sufficient to learn deep relu
  networks?
\newblock {\em CoRR}, arXiv:1911.12360, 2019.

\bibitem{CHERIDITO2021101540}
Patrick Cheridito, Arnulf Jentzen, and Florian Rossmannek.
\newblock Non-convergence of stochastic gradient descent in the training of
  deep neural networks.
\newblock {\em Journal of Complexity}, 64:101540, 2021.

\bibitem{ding2021overparameterization}
Zhiyan Ding, Shi Chen, Qin Li, and Stephen Wright.
\newblock Overparameterization of deep resnet: zero loss and mean-field
  analysis.
\newblock {\em arxiv:2105.14417}, 2021.

\bibitem{doi:10.1002/cnm.1640100303}
M.~W. M.~G. Dissanayake and N.~Phan-Thien.
\newblock Neural-network-based {A}pproximations for {S}olving {P}artial
  {D}ifferential {E}quations.
\newblock {\em Comm. Numer. Methods Engrg.}, 10:195--201, 1994.

\bibitem{Du2018}
S.~S. Du, X.~Zhai, B.~Poczos, and A.~Singh.
\newblock {Gradient descent provably optimizes over-parameterized neural
  networks}.
\newblock {\em arXiv e-prints}, arXiv:1810.02054, 2018.

\bibitem{duan2021convergence}
Chenguang Duan, Yuling Jiao, Yanming Lai, Xiliang Lu, and Zhijian Yang.
\newblock Convergence rate analysis for deep ritz method.
\newblock {\em arxiv:2103.13330}, 2021.

\bibitem{E2020}
Weinan E, Chao Ma, Stephan Wojtowytsch, and Lei Wu.
\newblock Towards a mathematical understanding of neural network-based machine
  learning: what we know and what we don't.
\newblock {\em arXiv e-prints}, arXiv:2009.10713, 2020.

\bibitem{Weinan2019}
Weinan E, Chao Ma, and Lei Wu.
\newblock A priori estimates of the population risk for two-layer neural
  networks.
\newblock {\em Communications in Mathematical Sciences}, 17(5):1407 -- 1425,
  2019.

\bibitem{frank2005nonlinear}
Till~Daniel Frank.
\newblock {\em Nonlinear Fokker-Planck equations: fundamentals and
  applications}.
\newblock Springer Science \& Business Media, 2005.

\bibitem{GUHRING2021107}
Ingo Gühring and Mones Raslan.
\newblock Approximation rates for neural networks with encodable weights in
  smoothness spaces.
\newblock {\em Neural Networks}, 134:107--130, 2021.

\bibitem{Gu2020}
Y.~Gu, H.~Yang, and C.~Zhou.
\newblock {SelectNet: Self-paced Learning for High-dimensional Partial
  Differential Equations}.
\newblock {\em Journal of Computational Physics}, 441:110444, 2021.

\bibitem{Han2018}
J.~Han, A.~Jentzen, and W.~E.
\newblock Solving high-dimensional partial differential equations using deep
  learning.
\newblock {\em Proc. Natl. Acad. Sci.}, 115(34):8505--8510, 2018.

\bibitem{he2016deep}
K.~{He}, X.~{Zhang}, S.~{Ren}, and J.~{Sun}.
\newblock Deep residual learning for image recognition.
\newblock In {\em 2016 IEEE Conference on Computer Vision and Pattern
  Recognition (CVPR)}, pages 770--778, 2016.

\bibitem{hess1976fokker}
Siegfried Hess.
\newblock Fokker-planck-equation approach to flow alignment in liquid crystals.
\newblock {\em Zeitschrift f{\"u}r Naturforschung A}, 31(9):1034--1037, 1976.

\bibitem{hon2021simultaneous}
Sean Hon and Haizhao Yang.
\newblock Simultaneous neural network approximations in sobolev spaces.
\newblock {\em arxiv:2109.00161}, 2021.

\bibitem{huggins2017quantifying}
Jonathan Huggins and James Zou.
\newblock Quantifying the accuracy of approximate diffusions and markov chains.
\newblock In {\em Artificial Intelligence and Statistics}, pages 382--391.
  PMLR, 2017.

\bibitem{HJKN19_814}
M.~Hutzenthaler, A.~Jentzen, Th. Kruse, and T.~A. Nguyen.
\newblock A proof that rectified deep neural networks overcome the curse of
  dimensionality in the numerical approximation of semilinear heat equations.
\newblock Technical Report 2019-10, Seminar for Applied Mathematics, ETH
  Z{\"u}rich, Switzerland, 2019.

\bibitem{hwang1994nonparametric}
Jenq-Neng Hwang, Shyh-Rong Lay, and Alan Lippman.
\newblock Nonparametric multivariate density estimation: a comparative study.
\newblock {\em IEEE Transactions on Signal Processing}, 42(10):2795--2810,
  1994.

\bibitem{iancu2001nonlinear}
Edmond Iancu, Andrei Leonidov, and Larry McLerran.
\newblock Nonlinear gluon evolution in the color glass condensate: I.
\newblock {\em Nuclear Physics A}, 692(3-4):583--645, 2001.

\bibitem{Ishimori1982}
Yuji Ishimori.
\newblock Solitons in a one-dimensional lennard-jones lattice.
\newblock {\em Progress of Theoretical Physics}, 68:402--410, 1982.

\bibitem{DBLP:journals/corr/abs-1806-07572}
Arthur Jacot, Franck Gabriel, and Cl{\'{e}}ment Hongler.
\newblock Neural tangent kernel: Convergence and generalization in neural
  networks.
\newblock {\em CoRR}, abs/1806.07572, 2018.

\bibitem{Jiao2021}
Y.~Jiao, G.~Shen, Y.~Lin, and J.~Huang.
\newblock Deep nonparametric regression on approximately low-dimensional
  manifolds.
\newblock {\em arXiv e-prints}, arXiv:2104.06708, 2021.

\bibitem{Sirignano2018}
S.~Justin and S.~Konstantinos.
\newblock Dgm: A deep learning algorithm for solving partial differential
  equations.
\newblock {\em J. Comput. Phys.}, 375:1339--1364, 2018.

\bibitem{Khoo2017SolvingPP}
Yuehaw Khoo, Jianfeng Lu, and Lexing Ying.
\newblock Solving parametric pde problems with artificial neural networks.
\newblock {\em arXiv: Numerical Analysis}, 2017.

\bibitem{KingmaB14}
D.~P. Kingma and J.~Ba.
\newblock {Adam: a method for stochastic optimization}.
\newblock {\em arXiv e-prints}, arXiv:1412.6980, 2014.

\bibitem{kumar2006solution}
Pankaj Kumar and S~Narayanan.
\newblock Solution of fokker-planck equation by finite element and finite
  difference methods for nonlinear systems.
\newblock {\em Sadhana}, 31(4):445--461, 2006.

\bibitem{712178}
I.E. Lagaris, A.~Likas, and D.~I. Fotiadis.
\newblock Artificial {N}eural {N}etworks for {S}olving {O}rdinary and {P}artial
  {D}ifferential {E}quations.
\newblock {\em IEEE Trans. Neural Networks}, 9:987--1000, 1998.

\bibitem{liang2021solving}
Senwei Liang, Shixiao~W. Jiang, John Harlim, and Haizhao Yang.
\newblock Solving pdes on unknown manifolds with machine learning.
\newblock {\em arxiv:2106.06682}, 2021.

\bibitem{liu2020neural}
Shu Liu, Wuchen Li, Hongyuan Zha, and Haomin Zhou.
\newblock Neural parametric fokker-planck equations.
\newblock {\em arXiv preprint arXiv:2002.11309}, 2020.

\bibitem{Liu2020Jul}
Ziqi Liu, Wei Cai, and Zhi-Qin~John Xu.
\newblock Multi-scale deep neural network (mscalednn) for solving
  poisson-boltzmann equation in complex domains.
\newblock {\em Communications in Computational Physics}, 28(5):1970–2001, Jun
  2020.

\bibitem{Shen3}
J.~Lu, Z.~Shen, H.~Yang, and S.~Zhang.
\newblock {Deep Network Approximation for Smooth Functions}.
\newblock {\em arXiv e-prints}, arXiv:2001.03040, 2020.

\bibitem{lu2021priori}
Jianfeng Lu, Yulong Lu, and Min Wang.
\newblock A priori generalization analysis of the deep ritz method for solving
  high dimensional elliptic equations.
\newblock {\em arxiv:2101.01708}, 2021.

\bibitem{bsp:2013}
Luo Lu, Hui Jiang, and Wing~H. Wong.
\newblock Multivariate density estimation by bayesian sequential partitioning.
\newblock {\em Journal of the American Statistical Association},
  108(504):1402--1410, 2013.

\bibitem{lu2020meanfield}
Yiping Lu, Chao Ma, Yulong Lu, Jianfeng Lu, and Lexing Ying.
\newblock A mean-field analysis of deep resnet and beyond: Towards provable
  optimization via overparameterization from depth.
\newblock {\em arxiv:2003.05508}, 2020.

\bibitem{LuoYang2020}
Tao Luo and Haizhao Yang.
\newblock {Two-Layer Neural Networks for Partial Differential Equations:
  Optimization and Generalization Theory}.
\newblock {\em arXiv e-prints}, arXiv:2006.15733, 2020.

\bibitem{mattingly2002ergodicity}
Jonathan~C Mattingly, Andrew~M Stuart, and Desmond~J Higham.
\newblock Ergodicity for {SDEs} and approximations: locally {Lipschitz} vector
  fields and degenerate noise.
\newblock {\em Stochastic processes and their applications}, 101(2):185--232,
  2002.

\bibitem{MeiE7665}
Song Mei, Andrea Montanari, and Phan-Minh Nguyen.
\newblock A mean field view of the landscape of two-layer neural networks.
\newblock {\em Proceedings of the National Academy of Sciences},
  115(33):E7665--E7671, 2018.

\bibitem{mishra2020estimates}
Siddhartha Mishra and Roberto Molinaro.
\newblock Estimates on the generalization error of physics informed neural
  networks (pinns) for approximating pdes.
\newblock {\em arxiv:2006.16144}, 2020.

\bibitem{misra2019mish}
Diganta Misra.
\newblock Mish: A self regularized non-monotonic activation function.
\newblock {\em arXiv preprint arXiv:1908.08681}, 2019.

\bibitem{Hadrien}
H.~Montanelli and Q.~Du.
\newblock New error bounds for deep networks using sparse grids.
\newblock {\em arXiv e-prints}, arXiv:1712.08688, 2017.

\bibitem{MO}
H.~Montanelli and H.~Yang.
\newblock Error bounds for deep relu networks using the kolmogorov–arnold
  superposition theorem.
\newblock {\em Neural Networks}, 129:1--6, 2020.

\bibitem{bandlimit}
Hadrien Montanelli, Haizhao Yang, and Qiang Du.
\newblock Deep relu networks overcome the curse of dimensionality for
  bandlimited functions.
\newblock {\em arXiv preprint arXiv:1903.00735}, 2019.

\bibitem{JMLR:v21:20-002}
Ryumei Nakada and Masaaki Imaizumi.
\newblock Adaptive approximation and generalization of deep neural network with
  intrinsic dimensionality.
\newblock {\em Journal of Machine Learning Research}, 21(174):1--38, 2020.

\bibitem{papamakarios2017masked}
George Papamakarios, Theo Pavlakou, and Iain Murray.
\newblock Masked autoregressive flow for density estimation.
\newblock In {\em Advances in Neural Information Processing Systems}, pages
  2338--2347, 2017.

\bibitem{poggio2017}
T.~Poggio, H.N. Mhaskar, L.~Rosasco, B.~Miranda, and Q.~Liao.
\newblock Why and when can deep---but not shallow---networks avoid the curse of
  dimensionality: {A} review.
\newblock {\em International Journal of Automation and Computing}, 14:503--519,
  2017.

\bibitem{RAISSI2019686}
M.~Raissi, P.~Perdikaris, and G.E. Karniadakis.
\newblock Physics-informed neural networks: A deep learning framework for
  solving forward and inverse problems involving nonlinear partial differential
  equations.
\newblock {\em Journal of Computational Physics}, 378:686 -- 707, 2019.

\bibitem{risken1996fokker}
Hannes Risken.
\newblock Fokker-planck equation.
\newblock In {\em The Fokker-Planck Equation}, pages 63--95. Springer, 1996.

\bibitem{rosenblatt1956remarks}
Murray Rosenblatt.
\newblock Remarks on some nonparametric estimates of a density function.
\newblock {\em The Annals of Mathematical Statistics}, pages 832--837, 1956.

\bibitem{10.1214/19-AOS1875}
Johannes Schmidt-Hieber.
\newblock {Nonparametric regression using deep neural networks with ReLU
  activation function}.
\newblock {\em The Annals of Statistics}, 48(4):1875 -- 1897, 2020.

\bibitem{sepehrian2015numerical}
Behnam Sepehrian and Marzieh~Karimi Radpoor.
\newblock Numerical solution of non-linear fokker--planck equation using finite
  differences method and the cubic spline functions.
\newblock {\em Applied mathematics and computation}, 262:187--190, 2015.

\bibitem{Shai2014}
Shai Shalev-Shwartz and Shai Ben-David.
\newblock {\em Understanding machine learning: From theory to algorithms}.
\newblock Cambridge university press, 2014.

\bibitem{Shen4}
Zuowei Shen, Haizhao Yang, and Shijun Zhang.
\newblock Deep network with approximation error being reciprocal of width to
  power of square root of depth.
\newblock {\em Neural Computation}, 2021.

\bibitem{siegel2021characterization}
Jonathan~W. Siegel and Jinchao Xu.
\newblock Characterization of the variation spaces corresponding to shallow
  neural networks.
\newblock {\em arxiv:2106.15002}, 2021.

\bibitem{siegel2021improved}
Jonathan~W. Siegel and Jinchao Xu.
\newblock Improved approximation properties of dictionaries and applications to
  neural networks.
\newblock {\em arxiv:2101.12365}, 2021.

\bibitem{spencer1993numerical}
BF~Spencer and LA~Bergman.
\newblock On the numerical solution of the fokker-planck equation for nonlinear
  stochastic systems.
\newblock {\em Nonlinear Dynamics}, 4(4):357--372, 1993.

\bibitem{tropp2012user}
Joel~A Tropp.
\newblock User-friendly tail bounds for sums of random matrices.
\newblock {\em Foundations of computational mathematics}, 12(4):389--434, 2012.

\bibitem{uria2016neural}
Benigno Uria, Marc-Alexandre C{\^o}t{\'e}, Karol Gregor, Iain Murray, and Hugo
  Larochelle.
\newblock Neural autoregressive distribution estimation.
\newblock {\em The Journal of Machine Learning Research}, 17(1):7184--7220,
  2016.

\bibitem{uy2020neural}
Wayne Isaac~T Uy and Mircea~D Grigoriu.
\newblock Neural network representation of the probability density function of
  diffusion processes.
\newblock {\em Chaos: An Interdisciplinary Journal of Nonlinear Science},
  30(9):093118, 2020.

\bibitem{wang2019nonparametric}
Zhipeng Wang and David~W Scott.
\newblock Nonparametric density estimation for high-dimensional
  data—algorithms and applications.
\newblock {\em Wiley Interdisciplinary Reviews: Computational Statistics},
  11(4):e1461, 2019.

\bibitem{xu2020solving}
Yong Xu, Hao Zhang, Yongge Li, Kuang Zhou, Qi~Liu, and J{\"u}rgen Kurths.
\newblock Solving fokker-planck equation using deep learning.
\newblock {\em Chaos: An Interdisciplinary Journal of Nonlinear Science},
  30(1):013133, 2020.

\bibitem{yarotsky2019}
Dmitry {Yarotsky} and Anton {Zhevnerchuk}.
\newblock The phase diagram of approximation rates for deep neural networks.
\newblock {\em arXiv e-prints}, page arXiv:1906.09477, June 2019.

\bibitem{Zhu2019}
Z.~Song Z.~A.-Zhu, Y.~Li.
\newblock A convergence theory for deep learning via over-parameterization.
\newblock In Kamalika Chaudhuri and Ruslan Salakhutdinov, editors, {\em
  Proceedings of the 36th International Conference on Machine Learning},
  volume~97 of {\em Proceedings of Machine Learning Research}, pages 242--252,
  Long Beach, California, USA, 2019. PMLR.

\bibitem{Zang2020}
Y.~Zang, G.~Bao, X.~Ye, and H.~Zhou.
\newblock Weak adversarial networks for high-dimensional partial differential
  equations.
\newblock {\em J. Comput. Phys.}, 411:109409, 2020.

\bibitem{zhai2020deep}
Jiayu Zhai, Matthew Dobson, and Yao Li.
\newblock A deep learning method for solving fokker-planck equations.
\newblock {\em arXiv preprint arXiv:2012.10696}, 2020.

\bibitem{zhl_fods:2020}
He~Zhang, John Harlim, and Xiantao Li.
\newblock Estimating linear response statistics using orthogonal polynomials:
  An rkhs formulation.
\newblock {\em Foundations of Data Science}, 2(4):443--485, 2020.

\bibitem{zhang2021error}
He~Zhang, John Harlim, and Xiantao Li.
\newblock Error bounds of the invariant statistics in machine learning of
  ergodic it\^o diffusions.
\newblock {\em Physica D (in press), arXiv preprint arXiv:2105.10102}, 2021.

\end{thebibliography}

\appendix
\section{Proofs for Section \ref{sec:convergence}}\label{appendixA}
\begin{proof}[Proof of Lemma \ref{lem02}]
Since $f_0\in\BRLU$, by \cite[Theorem 12]{E2020}, there exists a two-layer ReLU FNN $f^*$ with width $W$ such that $\|f^*\|_{L^\infty([0,1]^d)}\leq \|f_0\|_{\BRLU}$ and
\BEA
\|f^*-f_0\|_{L^\infty([0,1]^d)}\leq 4\|f_0\|_{\BRLU}(d+1)^{\frac{1}{2}}W^{-\frac{1}{2}}\leq4\sqrt{2}\|f_0\|_{\BRLU}d^{\frac{1}{2}}W^{-\frac{1}{2}}.\notag
\EEA
So $f^*\in\mathcal{F}_{2,W,\text{ReLU}}^P$. Since $\nu$ is absolutely continuous with respect to the Lebesgue measure, it follows that
\BEA\label{30}
\|f^*-f_0\|_{L^2_\nu([0,1]^d)}^2\leq 32\|f_0\|_{\BRLU}^2dW^{-1}.
\EEA

Also, \cite[Lemma 3.2]{Jiao2021} implies that
\BEA\label{31}
\mathbb{E}_\nu\left[|f_\text{NN}(\cdot,\Btheta^{f_0})-f_0|^2\right]\leq C\left[P^2W(d+W)\log(Wd+W^2)(\log N)^3N^{-1}
+\inf_{f\in\mathcal{F}_{2,W,\text{ReLU}}^P}\mathbb{E}_\nu\left[|f-f_0|^2\right]\right],
\EEA
where $C$ is a constant that does not depend on $d$, $N$, $W$, $f_0$, $P$. Combining \eqref{30} and \eqref{31} completes the proof.
\end{proof}

\begin{proof}[Proof of Lemma \ref{thm02}]
Denote $\hat{e}:=q-\hat{p}$. On one hand, using integration by parts,
\BEA\label{11}
\int_\Omega\hat{\mathcal{L}}^*\hat{e}\cdot\hat{e}\,\td \Bx&\geq&\int_\Omega\sum_{i,j=1}^d\frac{1}{2}B_\tNN^{ij}\hat{e}_{x_i}\hat{e}_{x_j}\td \Bx-\int_{\partial\Omega}\left(\sum_{i,j=1}^d\frac{1}{2}B_\tNN^{ij}|\hat{e}_{x_i}|\cdot|\tn_j|\right)|\hat{e}|\td s
+\int_\Omega\sum_{i=1}^da_\tNN^i\hat{e}_{x_i}\cdot\hat{e}
+\left(\sum_{i=1}^d\frac{\partial a_\tNN^i}{\partial x_i}\right)\hat{e}^2\td \Bx\notag \\
&\geq&\frac{1}{2}\Lambda\int_\Omega\|\nabla\hat{e}\|^2\td x-\frac{1}{2}dB_1\int_{\partial\Omega}\|\nabla\hat{e}\|\cdot|\hat{e}|\td s \notag \\
&\geq&\frac{1}{2}\Lambda\|\nabla\hat{e}\|_{L^2(\Omega)}^2-\frac{1}{2}dB_1\left(\|\nabla q\|_{L^2(\partial\Omega)}+\|\nabla\hat{p}\|_{L^2(\partial\Omega)}\right)\|\hat{e}\|_{L^2(\partial\Omega)}\notag\\
&\geq&\frac{1}{2}\Lambda\|\nabla\hat{e}\|_{L^2(\Omega)}^2-\frac{1}{2}dB_1\left(B_2+\epsilon_{\hat{p}}\right)\|\hat{e}\|_{L^2(\partial\Omega)},
\EEA
where $\tn_j$ is the $j$-th component of the outward unit normal vector.

On the other hand,
\begin{equation}\label{12}
\int_\Omega\hat{\mathcal{L}}^*\hat{e}\cdot\hat{e}\td x=\int_\Omega\hat{\mathcal{L}}^*q\cdot\hat{e}\td x\leq\|\hat{\mathcal{L}}^*q\|_{L^2(\Omega)}\cdot\|\hat{e}\|_{L^2(\Omega)}.
\end{equation}
Combining \eqref{11} and \eqref{12} leads to
\begin{equation}\label{13}
\|\nabla\hat{e}\|_{L^2(\Omega)}^2\leq2\Lambda^{-1}\|\hat{\mathcal{L}}^*q\|_{L^2(\Omega)}\cdot\|\hat{e}\|_{L^2(\Omega)}+\Lambda^{-1}dB_1(B_2+\epsilon_{\hat{p}})\|\hat{e}\|_{L^2(\partial\Omega)}.
\end{equation}

Next, by Poincar\'{e} inequality, there exists some $C_1>0$ that only depends on $\Omega$ such that
\begin{equation*}
\left\|\hat{e}-|\Omega|^{-1}\int_\Omega\hat{e}\td x\right\|_{L^2(\Omega)}\leq C_1\|\nabla\hat{e}\|_{L^2(\Omega)},
\end{equation*}
which leads to
\begin{equation*}
\|\hat{e}\|_{L^2(\Omega)}\leq |\Omega|^{-1}\left|\int_\Omega\hat{e}\td x \right|\|1\|_{L^2(\Omega)}+ C_1\|\nabla\hat{e}\|_{L^2(\Omega)}
\leq C_2\left(\left|\int_\Omega\hat{e}\td x\right|+\|\nabla\hat{e}\|_{L^2(\Omega)}\right),
\end{equation*}
where $C_2=\max(C_1,|\Omega|^{-1/2})$. Therefore, by \eqref{13} and the fact $\int_\Omega\hat{p}\td \Bx=1$,
\BEA\label{14}
\|\hat{e}\|_{L^2(\Omega)}^2 &\leq& C_3\left[\left|\int_\Omega\hat{e}\td \Bx\right|^2+\|\nabla\hat{e}\|_{L^2(\Omega)}^2\right]\notag\\
&\leq& C_3\left[\left|\int_\Omega q\td x-1\right|^2+2\Lambda^{-1}\|\hat{\mathcal{L}}^*q\|_{L^2(\Omega)}\cdot\|\hat{e}\|_{L^2(\Omega)}+\Lambda^{-1}dB_1(B_2+\epsilon_{\hat{p}})\|\hat{e}\|_{L^2(\partial\Omega)}\right],
\EEA
where $C_3=2C_2^2$. Using the Young's inequality $2C_3\Lambda^{-1}\|\hat{\mathcal{L}}^*q\|_{L^2(\Omega)}\cdot\|\hat{e}\|_{L^2(\Omega)}\leq\frac{4C_3^2\Lambda^{-2}\|\hat{\mathcal{L}}^*q\|_{L^2(\Omega)}^2+\|\hat{e}\|_{L^2(\Omega)}^2}{2}$, it follows from \eqref{14} that
\begin{equation}\label{15}
\frac{1}{2}\|\hat{e}\|_{L^2(\Omega)}^2
\leq C_3\left|\int_\Omega q\td x-1\right|^2+
2C_3^2\Lambda^{-2}\|\hat{\mathcal{L}}^*q\|_{L^2(\Omega)}^2+C_3\Lambda^{-1}dB_1(B_2+\epsilon_{\hat{p}})\|\hat{e}\|_{L^2(\partial\Omega)}.
\end{equation}

Note $\|\hat{e}\|_{L^2(\partial\Omega)}\leq\|\hat{p}\|_{L^2(\partial\Omega)}+\|q\|_{L^2(\partial\Omega)}\leq\epsilon_{\hat{p}}+\|q\|_{L^2(\partial\Omega)}$, it follows from \eqref{15} that
\BEA
\|\hat{e}\|_{L^2(\Omega)}^2
&\leq& 2C_3\left|\int_\Omega q\td x-1\right|^2
+4C_3^2\Lambda^{-2}\|\hat{\mathcal{L}}^*q\|_{L^2(\Omega)}^2
+2C_3\Lambda^{-1}dB_1(B_2+\epsilon_{\hat{p}})\epsilon_{\hat{p}} \notag \\ && +2C_3\Lambda^{-1}dB_1(B_2+\epsilon_{\hat{p}})\|q\|_{L^2(\partial\Omega)} \notag \\
&\leq& C\left(J[q]+d(B_2+\epsilon_{\hat{p}})J[q]^\frac{1}{2}+d(B_2+\epsilon_{\hat{p}})\epsilon_{\hat{p}}\right),\notag
\EEA
where $C$ only depends on $\Omega,\Lambda,B_1,\lambda_1,\lambda_2$.
\end{proof}

\begin{proof}[Proof of Lemma \ref{lem03}]
Let $f=\mathbb{E}_{(c,\Bw)\sim\rho}[c\dsigma(\Bw^\top\Bx)]$ for some $\rho$ taking the infimum in \eqref{07}. Then $\hat{\OL}^*f=\mathbb{E}_{(c,\Bw)\sim\rho}[\hat{\OL}^*(c\dsigma(\Bw^\top\Bx))]$. Using the homogeneity of the neuron $c\dsigma(\Bw^\top\Bx)$, we may assume that $\|\Bw\|_1=1$ and $|c|=\|f\|_{\Bdsg}$ $\rho$-almost everywhere. Indeed, denote $p_0$ as the density of $\rho$, we define the probability measure $\rho^*$ with the density
\BEA
p_0^*(\hat{c},\hat{\Bw})=\begin{cases}\int_{c\|\Bw\|_1^3=\hat{c}}p_0(c,\Bw)\td c\td \Bw,\quad\text{if}~\|\hat{\Bw}\|_1=1,\\0,\quad\text{otherwise},\end{cases}
\EEA
then it can be verified that $\rho^*\in P_f$, $\mathbb{E}_\rho|c|\|\Bw\|_1^3=\mathbb{E}_{\rho^*}|\hat{c}|\|\hat{\Bw}\|_1^3$ and $\text{supp}(p_0^*)\subset\mathbb{R}\times\{\|\hat{\Bw}\|_1=1\}$. Moreover, we define the probability measure $\rho^{**}$ with the density
\BEA
p_0^{**}(\tilde{c},\tilde{\Bw})=\begin{cases}\|f\|_{\Bdsg}^{-1}\int_0^{+\infty}|\hat{c}|p_0^*(\hat{c},\hat{\Bw})\td\hat{c},\quad\text{if}~\tilde{c}=\|f\|_{\Bdsg},\|\tilde{\Bw}\|_1=1,\\\|f\|_{\Bdsg}^{-1}\int_{-\infty}^0|\hat{c}|p_0^*(\hat{c},\hat{\Bw})\td\hat{c},\quad\text{if}~\tilde{c}=-\|f\|_{\Bdsg},\|\tilde{\Bw}\|_1=1,\\0,\quad\text{otherwise},\end{cases}
\EEA
then it can be verified that $\rho^{**}\in P_f$, $\mathbb{E}_{\rho^*}|\hat{c}|\|\hat{\Bw}\|_1^3=\mathbb{E}_{\rho^{**}}|\tilde{c}|\|\tilde{\Bw}\|_1^3$ and $\text{supp}(p_0^{**})\subset\{\tilde{c}=\pm\|f\|_{\Bdsg}\}\times\{\|\tilde{\Bw}\|_1=1\}$.

Let $\{(c_m,\Bw_m)\}$ be $M$ independent and identically distributed samples with $\rho$. By \cite[Lemma 26.2]{Shai2014},
\begin{multline}\label{09}
\mathbb{E}_{\{(c_m,\Bw_m)\}\sim\rho^M}\left[\underset{\Bx\in\Omega}{\sup}~\hat{\OL}^*\left(\frac{1}{M}\sum_{m=1}^Mc_m\dsigma(\Bw_m^\top\Bx)\right)-\hat{\OL}^*f(\Bx)\right]\\
=\mathbb{E}_{\{(c_m,\Bw_m)\}\sim\rho^M}\left[\underset{\Bx\in\Omega}{\sup}\left(\frac{1}{M}\sum_{m=1}^M~\hat{\OL}^*c_m\dsigma(\Bw_m^\top\Bx)-\mathbb{E}_{(c,\Bw)\sim\rho}[\hat{\OL}^*(c\dsigma(\Bw^\top\Bx))]\right)\right]\\ \leq
2\mathbb{E}_{\{(c_m,\Bw_m)\}\sim\rho^M}\mathbb{E}_\tau\left[\underset{\Bx\in\Omega}{\sup}~\frac{1}{M}\sum_{m=1}^M\tau_m\hat{\OL}^*(c_m\dsigma(\Bw_m^\top\Bx))\right],
\end{multline}
where $\tau_m=\pm1$ with probability $1/2$ are independent Rademacher variables.

Note that
\begin{multline}\label{08}
 \mathbb{E}_\tau\left[\underset{\Bx\in\Omega}{\sup}~\frac{1}{M}\sum_{m=1}^M\tau_m\hat{\OL}^*(c_m\dsigma(\Bw_m^\top\Bx))\right]\\
=\mathbb{E}_\tau\left[\underset{\Bx\in\Omega}{\sup}~\frac{1}{M}\sum_{m=1}^M\tau_mc_m\left(\frac{1}{2}\Bw_m^\top\BB_\tNN\Bw_m \dsigma''(\Bw_m^\top\Bx)+\Ba_\tNN^\top\Bw_m\dsigma'(\Bw_m^\top\Bx)+\left(\sum_{i=1}^d\frac{\partial a_\tNN^i}{\partial x_i}\right)\dsigma(\Bw_m^\top\Bx)
\right)\right] \\
\leq\mathbb{E}_\tau\left[\underset{\Bx\in\Omega}{\sup}~\frac{1}{M}\sum_{m=1}^M\frac{1}{2}\tau_mc_m\Bw_m^\top\BB_\tNN\Bw_m \dsigma''(\Bw_m^\top\Bx)\right]
+\mathbb{E}_\tau\left[\underset{\Bx\in\Omega}{\sup}~\frac{1}{M}\sum_{m=1}^M\tau_mc_m\Ba_\tNN^\top\Bw_m\dsigma'(\Bw_m^\top\Bx)\right]\hspace{10pt}\\
+\mathbb{E}_\tau\left[\underset{\Bx\in\Omega}{\sup}~\frac{1}{M}\sum_{m=1}^M\tau_mc_m\left(\sum_{i=1}^d\frac{\partial a_\tNN^i}{\partial x_i}\right)\dsigma(\Bw_m^\top\Bx)\right]\hspace{2in}
\end{multline}

For the first term in \eqref{08}, by the contraction lemma for Rademacher complexities \cite[Lemma 26.9]{Shai2014}, we have
\begin{multline}\label{16}
\mathbb{E}_\tau\left[\underset{\Bx\in\Omega}{\sup}~\frac{1}{M}\sum_{m=1}^M\frac{1}{2}\tau_mc_m\Bw_m^\top\BB_\tNN\Bw_m \dsigma''(\Bw_m^\top\Bx)\right]
=\mathbb{E}_\tau\left[\underset{\Bx\in\Omega}{\sup}~\frac{1}{M}\sum_{m=1}^M\tau_m\dsigma''(\frac{1}{2}c_m\Bw_m^\top\BB_\tNN\Bw_m\cdot\Bw_m^\top\Bx)\right]\\
\leq\mathbb{E}_\tau\left[\underset{\Bx\in\Omega}{\sup}~\frac{1}{M}\sum_{m=1}^M\tau_mB_1|c_m|\|\Bw_m\|_1^2\cdot\Bw_m^\top\Bx\right]
=\frac{B_1}{M}\mathbb{E}_\tau\left[\underset{\Bx\in\Omega}{\sup}~\Bx^\top\sum_{m=1}^M\tau_m|c_m|\|\Bw_m\|_1^2\cdot\Bw_m\right]\\
\leq B_1\mathbb{E}_\tau\left\|\frac{1}{M}\sum_{m=1}^M\tau_m|c_m|\|\Bw_m\|_1^2\cdot\Bw_m\right\|_1.
\end{multline}

Similarly, we can derive
\BEA\label{17}
\mathbb{E}_\tau\left[\underset{\Bx\in\Omega}{\sup}~\frac{1}{M}\sum_{m=1}^M\tau_mc_m\Ba_\tNN^\top\Bw_m\dsigma'(\Bw_m^\top\Bx)\right]
\leq B_1\mathbb{E}_\tau\left\|\frac{1}{2M}\sum_{m=1}^M\tau_m|c_m|\|\Bw_m\|_1\cdot\Bw_m\right\|_1
\EEA
and
\BEA\label{18}
\mathbb{E}_\tau\left[\underset{\Bx\in\Omega}{\sup}~\frac{1}{M}\sum_{m=1}^M\tau_mc_m\left(\sum_{i=1}^d\frac{\partial a_\tNN^i}{\partial x_i}\right)\dsigma(\Bw_m^\top\Bx)\right]
\leq B_1\mathbb{E}_\tau\left\|\frac{1}{2M}\sum_{m=1}^M\tau_m|c_m|\Bw_m\right\|_1.
\EEA
Denote $\Bu_m:=c_m\Bw_m$, then $\|\Bu_m\|_1=\|f\|_{\Bdsg}$. We combine \eqref{09}-\eqref{18} and obtain
\BEA
\mathbb{E}_{\{(c_m,\Bw_m)\}\sim\rho^M}\left[\underset{\Bx\in\Omega}{\sup}~\hat{\OL}\left(\frac{1}{M}\sum_{m=1}^Mc_m\dsigma(\Bw_m^\top\Bx)\right)-\hat{\OL}f(\Bx)\right]
&\leq& 2\underset{\|\Bu_m\|_1\leq\|f\|_{\Bdsg}}{\sup}~2B_1\mathbb{E}_\tau\|\frac{1}{M}\sum_{m=1}^M\tau_m\Bu_m\|_1 \notag\\
&\leq& 2\|f\|_{\Bdsg}\underset{\|\Bu_m\|_1\leq1}{\sup}~2B_1\mathbb{E}_\tau\|\frac{1}{M}\sum_{m=1}^M\tau_m\Bu_m\|_1\notag\\
&\leq& 2\sqrt{d}\|f\|_{\Bdsg}\underset{\|\Bu_m\|_2\leq1}{\sup}~2B_1\mathbb{E}_\tau\|\frac{1}{M}\sum_{m=1}^M\tau_m\Bu_m\|_2 \notag\\
&\leq& 4B_1\|f\|_{\Bdsg}\sqrt{d/M}
\EEA
by using the Rademacher complexity of the unit ball \cite[Lemma 26.10]{Shai2014}. Applying the same argument to $-\left(\hat{\OL}\left(\frac{1}{M}\sum_{m=1}^Mc_m\dsigma(\Bw_m^\top\Bx)\right)-\hat{\OL}f(\Bx)\right)$ leads to
\BEA
\mathbb{E}_{\{(c_m,\Bw_m)\}\sim\rho^M}\left[\underset{\Bx\in\Omega}{\sup}~\left|\hat{\OL}\left(\frac{1}{M}\sum_{m=1}^Mc_m\dsigma(\Bw_m^\top\Bx)\right)-\hat{\OL}f(\Bx)\right|\right]
\leq4B_1\|f\|_{\Bdsg}\sqrt{d/M}.
\EEA

By similar argument, we can derive
\BEA
\mathbb{E}_{(c,\Bw)\sim\rho}\left[\underset{\Bx\in\Omega}{\sup}\left|\frac{1}{M}\sum_{m=1}^Mc_m\dsigma(\Bw_m^\top\Bx)-f(\Bx)\right|\right],~~\mathbb{E}_{(c,\Bw)\sim\rho}\left[\underset{\Bx\in\partial\Omega}{\sup}\left|\frac{1}{M}\sum_{m=1}^Mc_m\dsigma(\Bw_m^\top\Bx)-f(\Bx)\right|\right]\leq\|f\|_{\Bdsg}\sqrt{d/M}.
\EEA

Therefore we have
\begin{multline}
\mathbb{E}_{(c,\Bw)\sim\rho^M}\Bigg[\underset{\Bx\in\Omega}{\sup}~\left|\hat{\OL}\left(\frac{1}{M}\sum_{m=1}^Mc_m\dsigma(\Bw_m^\top\Bx)\right)-\hat{\OL}f(\Bx)\right|
+\underset{\Bx\in\Omega}{\sup}~\left|\frac{1}{M}\sum_{m=1}^Mc_m\dsigma(\Bw_m^\top\Bx)-f(\Bx)\right|\\
+\underset{\Bx\in\partial\Omega}{\sup}~\left|\frac{1}{M}\sum_{m=1}^Mc_m\dsigma(\Bw_m^\top\Bx)-f(\Bx)\right|\Bigg]
\leq \left(4B_1+2\right)\|f\|_{\Bdsg}\sqrt{d/M},
\end{multline}
which implies there exists $\{(c_m,\Bw_m)\}_{m=1}^M$ such that the inequality holds. Then the FNN $\sum_{m=1}^M(c_m/M)\dsigma(\Bw_m^\top\Bx)\in\mathcal{F}_{2,M,\dsigma,\max\left\{\|f\|_{\Bdsg}/M,1\right\}}$ satisfies \eqref{10}.
\end{proof}

\begin{proof}[Proof of Lemma \ref{thm03}]
Denote $\hpNN^S(\Bx)=\hpNN(\Bx;\Btheta^S)$. Since $\hpNN\in\mathcal{F}_{2,M,\dsigma,Q}$, using the expression in \eqref{two_layer_FNN} we have $\|\nabla\hpNN^S\|_{L^2(\partial\Omega)}\leq\frac{1}{2}MQ^4|\partial\Omega|^\frac{1}{2}=\frac{1}{2}MQ^4(2d)^\frac{1}{2}$. Then the inequality \eqref{25} directly follows Lemma \ref{thm02}. For the rest, we use $C$ to represent any constant which on depends on $\Omega$, $\Lambda$, $B_1$, $\lambda_1$ and $\lambda_2$. On one hand,
\BEA\label{23}
|J[\hpNN^S]-J_S[\hpNN^S]|&\leq&\left|\|\OL\hpNN^S\|_{L^2(\Omega)}^2-\frac{|\Omega|}{N_1}\sum_{n=1}^{N_1}|\OL\hpNN^S(\Bx_{\rm I}^n)|^2\right|
+\lambda_2\left|\|\hpNN^S\|_{L^2(\partial\Omega)}^2-\frac{|\partial\Omega|}{N_3}\sum_{n=1}^{N_3}|\hpNN^S(\Bx_{\rm III}^n)|^2\right|\notag\\
&&+\lambda_1\left|\int_\Omega\hpNN^S(\Bx)\td\Bx-\frac{|\Omega|}{N_2}\sum_{n=1}^{N_2}\hpNN^S(\Bx_{\rm II}^n)\right|\cdot\left|\int_\Omega\hpNN^S(\Bx)\td\Bx+\frac{|\Omega|}{N_2}\sum_{n=1}^{N_2}\hpNN^S(\Bx_{\rm II}^n)-2\right|.
\EEA
By virtue of \cite[Theorem 3.2]{LuoYang2020}, with probability at least $1-\delta/3$,
\BEA\label{24}
\left|\|\OL\hpNN^S\|_{L^2(\Omega)}^2-\frac{|\Omega|}{N_1}\sum_{n=1}^{N_1}|\OL\hpNN^S(\Bx_{\rm I}^n)|^2\right|
+\lambda_2\left|\|\hpNN^S\|_{L^2(\partial\Omega)}^2-\frac{|\partial\Omega|}{N_3}\sum_{n=1}^{N_3}|\hpNN^S(\Bx_{\rm III}^n)|^2\right|\leq CI_1.
\EEA

Similarly, by the fact $|\hpNN^S(\Bx)|\leq MQ^4/6$ for all $\Bx$ and Lemma \ref{lem04}, we have with probability at least $1-\delta/3$,
\BEA
\left|\int_\Omega\hpNN^S(\Bx)\td\Bx-\frac{|\Omega|}{N_2}\sum_{n=1}^{N_2}\hpNN^S(\Bx_{\rm II}^n)\right|\leq C MQ^4\sqrt{\log(6/\delta)/N_2},
\EEA
and $\frac{|\Omega|}{N_2}\sum_{n=1}^{N_2}\hpNN^S(\Bx_{\rm II}^n)\leq C MQ^4$. Then we have
\BEA\label{26}
\lambda_1\left|\int_\Omega\hpNN^S(\Bx)\td\Bx-\frac{|\Omega|}{N_2}\sum_{n=1}^{N_2}\hpNN^S(\Bx_{\rm II}^n)\right|\cdot\left|\int_\Omega\hpNN^S(\Bx)\td\Bx+\frac{|\Omega|}{N_2}\sum_{n=1}^{N_2}\hpNN^S(\Bx_{\rm II}^n)-2\right|\leq C I_2.
\EEA

On the other hand, by Lemma \ref{lem03} there exists some $p_\tNN\in F_{2,M,\dsigma,Q}$ such that
\BEA\label{20}
\underset{\Bx\in\Omega}{\sup}~\left|\hat{\OL}p_\tNN(\Bx)\right|
+\underset{\Bx\in\Omega}{\sup}~\left|p_\tNN(\Bx)-\hat{p}(\Bx)\right|
+\underset{\Bx\in\partial\Omega}{\sup}~\left|p_\tNN(\Bx)-\hat{p}(\Bx)\right|
\leq C\|\hat{p}\|_{\Bdsg}\sqrt{d/M}.
\EEA
Note that $\int_\Omega\hat{p}\td\Bx=1$, we have
\BEA\label{21}
\left|\frac{|\Omega|}{N_2}\sum_{n=1}^{N_2} p_\tNN(\Bx_{\rm II}^n)-1\right|^2\leq 2\left(\left|\frac{|\Omega|}{N_2}\sum_{n=1}^{N_2} \left(p_\tNN(\Bx_{\rm II}^n)-\hat{p}(\Bx_{\rm II}^n)\right)\right|^2+\left|\frac{|\Omega|}{N_2}\sum_{n=1}^{N_2} \hat{p}(\Bx_{\rm II}^n)-\int_\Omega\hat{p}\td\Bx\right|^2\right).
\EEA
and
\BEA\label{22}
\left|p_\tNN(\Bx_{\rm III}^n)\right|^2\leq2\left(\left|p_\tNN(\Bx_{\rm III}^n)-\hat{p}(\Bx_{\rm III}^n)\right|^2+\epsilon_{\hat{p}}^2\right),
\EEA
using the fact that $|\hat{p}(\Bx)|\leq \epsilon_{\hat{p}}$ on $\partial\Omega$ in Assumption~\ref{assum_03}.

Then it follows \eqref{20}-\eqref{22} and Lemma \ref{lem04} that with probability at least $1-\delta/3$
\BEA\label{27}
J_S[\hpNN^S]\leq J_S[p_\tNN]&\leq & \frac{1}{N_1}\sum_{n=1}^{N_1} \left|\OL p_\tNN(\Bx_{\rm I}^n)\right|^2
+2\lambda_1 \left|\frac{|\Omega|}{N_2}\sum_{n=1}^{N_2} \left(p_\tNN(\Bx_{\rm II}^n)-\hat{p}(\Bx_{\rm II}^n)\right)\right|^2
+2\lambda_1 \left|\frac{|\Omega|}{N_2}\sum_{n=1}^{N_2} \hat{p}(\Bx_{\rm II}^n)-\int_\Omega\hat{p}\right|^2 \notag\\
&&+2\lambda_2 \frac{|\partial\Omega|}{N_3}\sum_{n=1}^{N_3}\left(\left|p_\tNN(\Bx_{\rm III}^n)-\hat{p}(\Bx_{\rm III}^n)\right|^2+\epsilon_{\hat{p}}^2\right)
\leq CI_3.
\EEA

Finally, the proof can be completed by using \eqref{23}, \eqref{24}, \eqref{26}, \eqref{27} and the fact $J[\hpNN^S]\leq\left|J[\hpNN^S]-J_S[\hpNN^S]\right|+J_S[\hpNN^S]$.
\end{proof}

\end{document}